\documentclass{article}

\usepackage{amssymb,amsmath,mathrsfs}
\usepackage[colorlinks=false]{hyperref}
\usepackage{comment}
\usepackage{diagbox}

\usepackage{geometry}
\usepackage{diagbox}
\usepackage{manyfoot}
\usepackage{graphicx,color}
\usepackage[outercaption]{sidecap}
\usepackage{xfrac}
\usepackage{subfig}
\usepackage{tikz}
\usepackage{enumitem}
\usepackage{pdfsync}
\usepackage{fancyhdr}
\usepackage[nottoc,notlof,notlot]{tocbibind}
\usepackage[titles,subfigure]{tocloft}

\usepackage{amssymb}
\usepackage{amsmath}
\usepackage{latexsym}
\usepackage{amsthm}
\usepackage{eucal}
\usepackage{amsthm}
\theoremstyle{plain}
\newtheorem{theorem}{Theorem}[section]

\newtheorem{lemma}[theorem]{Lemma}
\newtheorem{proposition}[theorem]{Proposition}
\theoremstyle{definition}

\newtheorem{remark}[theorem]{Remark}

\DeclareMathOperator*{\R}{\mathbb{R}}

\DeclareMathOperator*{\dist}{dist}
\DeclareMathOperator*{\supp}{supp}
\DeclareMathOperator*{\loc}{loc}
\DeclareMathOperator*{\capi}{cap}
\DeclareMathOperator*{\reg}{reg}

\def\XXint#1#2#3{{\setbox0=\hbox{$#1{#2#3}{\int}$} 
  \vcenter{\hbox{$#2#3$}}\kern-.5\wd0}}

\numberwithin{equation}{section}

\title{Variational analysis of discrete Dirichlet problems in periodically perforated domains }
\author{
{\sc Giuliana Fusco}
\\ \small Scuola Superiore Meridionale,\\
 \small
via Mezzocannone 4, 80134 Napoli, Italy\\
\small g.fusco@ssmeridionale.it
}
\date{
}                                      

\begin{document}
\maketitle

\noindent
{\bf Abstract.}
In this paper we study the asymptotic behavior of a family of discrete functionals as the lattice size, $\varepsilon>0$, tends to zero. We consider pairwise interaction energies satisfying $p$-growth conditions, $p<d$, $d$ being the dimension of the reference configuration, defined on discrete functions subject to Dirichlet conditions on a $\delta$-periodic array of small squares of side $r_{\delta}\sim \delta^{d/d-p}$. Our analysis is performed in the framework of $\Gamma$-convergence and we prove that, in the regime $\varepsilon=o(r_{\delta})$, the discrete energy and their continuum counterpart share the same $\Gamma$-limit and the effect of the constraints leads to a capacitary term in the limit energy as in the classical theory of periodically perforated domains for local integral functionals.
\bigskip

$\Gamma$-convergence.
\smallskip


\section{Introduction}
In the last decades many efforts have been done to investigate the relation between continuum theories and general atomic interaction energies. The reason is related to their relevance for applications in different directions. On one hand the derivation of continuum models from atomistic systems provides a microscopical theoretical justification of continuum theories. On the other hand continuum models, obtained as an approximation of atomic interaction energies usually correspond to a coarse-graining description of the latters, which smooth away fine details but are able to capture the main features of the original problem, bringing a lot of advantages from a numerical view point. In this paper we are interested in studying the asymptotic behaviour of a specific class of atomic interaction energies satisfying suitable $p$-growth conditions, with $p<d$, $d$ being the dimension of the reference configuration, and  subject to a Dirichlet constraint on periodically perforated domains. In the continuous setting, the problem is modelled, in the simple case, by imposing a  Dirichlet condition on a periodic array $U_{\delta,R}$ of small balls of radius $R$ and centers on a $\delta$-periodic lattice. For instance, for a suitable choice  of the size $R=R_\delta$ of the perforations  minimum problems of the form 
\begin{equation*}
    \min\left\{\int_{\Omega}(|\nabla u|^{p}-gu)\, dx : u =0\  \text{on}\ U_{\delta,R}\right\}
\end{equation*}
can be approximated by
\begin{equation*}
    \min\left\{\int_{\Omega}(|\nabla u|^{p}+C_{p}|u|^{p}-gu)\, dx\right\},
\end{equation*}
where $C_{p}$ is the $p$-capacity of the unit ball $B_{1}(0)$ in $\mathbb{R}^{d}$ and the term $C_{p}|u|^{p}$ accounts for the energetic contribution near the perforations. This "extra term" appeared in literature for the first time in the well-known paper by D. Cioranescu and F. Murat \cite{CM}, where they study the scalar case for $p=2$. It must be noted that in the subcritical case $p<d$ the only relevant scaling is 
\begin{equation*}
    R\sim \delta^{\frac{d}{d-p}},
\end{equation*}
 since the other cases, i.e. $R=o\big(\delta^{\frac{d}{d-p}}\big)$ and $R=O\big(\delta^{\frac{d}{d-p}}\big)$, lead to trivial limit energy. In \cite{AB} N. Ansini and A. Braides set the problem in the framework of $\Gamma$-convergence and extended the analysis to a general class of energies defined on vector valued functions. Under suitable assumptions on the density function $f$, they prove that, as $\delta \rightarrow 0$ and $R_{\delta}=\delta^{d/d-p}$ functionals of the form
\begin{equation*}
    \int_{\Omega}f(\nabla u)\, dx,\ \quad u\in W^{1,p}(\Omega,\mathbb{R}^m),\  u=0\ \text{on}\ U_{\delta,R}
\end{equation*}
$\Gamma$-converge to the functional defined by
\begin{equation}
    \int_{\Omega}f(\nabla u)\, dx+\int_{\Omega}\varphi(u)\, dx,\ \quad u\in W^{1,p}(\Omega,\mathbb{R}^m),
    \label{2}
\end{equation}
where $\varphi(z)$ is described by the non linear capacitary formula
\begin{equation*}
 \varphi(z)=\inf\bigg\{\int_{\mathbb{R}^{d}}f(\nabla v)\, dx: v-z \in W^{1,p}_{\loc}(\mathbb{R}^{d};\mathbb{R}^{m})\cap L^{p^{*}}(\mathbb{R}^{d};\mathbb{R}^{m}), v=0\ \text{on}\ B_{1}(0)\bigg\}.
\end{equation*} 
In this paper we consider pairwise interaction energies of the form 
\begin{equation*}
F_{\varepsilon}(u)=\sum_{\xi \in \mathbb{Z}^{d}}\sum_{\alpha,\
\alpha+\varepsilon\xi\ \in\ \varepsilon\ \mathbb{Z}^{d}\cap\Omega}\varepsilon^{d}f\left(\xi,\frac{u(\alpha +\varepsilon\xi)-u(\alpha)}{\varepsilon|\xi|}\right)
\end{equation*}
defined on functions $u:\alpha \in \varepsilon\mathbb{Z}^{d}\cap \Omega \mapsto u(\alpha) \in \mathbb{R}^{m}$, identified with their piecewise constant interpolation on each cell of the lattice $\varepsilon\mathbb{Z}^{d}$. For any $\xi \in  \mathbb{Z}^{d}$, the energy densities $f(\xi,\cdot)$ are $p$-homogeneous, locally Lipschitz continuous, and satisfy suitable decay assumptions as $|\xi|\rightarrow +\infty$ (see hypotheses (H), (G), and (L) in Section $3$). Under this assumptions, as a particular case of the asymptotic analysis provided by R. Alicandro and M. Cicalese in \cite{alicic} (see Theorem \ref{51}), the $\Gamma$-limit of $F_{\varepsilon}$ is given by the local functional
\begin{equation*}
    \label{local functional}
    F_{0}(u)=\int_{\Omega}f_{\hom}(\nabla u)\, dx, \ \quad u\in W^{1,p}(\Omega;\mathbb{R}^m),
\end{equation*}
where $f_{\hom}(M)$ is given by a suitable homogenization formula. We then impose the admissible functions $u$ to satisfy the constraint $u=0$ on a periodic array $P_{\delta}$ of small squares of side $r_{\delta}\sim \delta^{d/d-p}$ and centers on a $\delta$-periodic lattice. We are interested in the asymptotic behaviour of $F_{\varepsilon}$ as $\varepsilon$ and $\delta$ go to $0$ according to the interplay between the scales $\varepsilon, \delta$. The case in which $\varepsilon \sim r_{\delta}$ has been studied by L. Sigalotti in \cite{Siga}, where she proves that limit energy is still of the form \eqref{2} where $f$ is replaced by $f_{\hom}$ and $\varphi$ is given by a non local capacitary formula. In this paper we focus our attention on the case $\varepsilon=o(r_{\delta})$ and we address the question whether or not $F_\varepsilon$ and their continuous approximation $F_0$ share the same asymptotic behaviour under imposed constraints on the admissible functions. We show that this the case. Indeed, in the main result of the paper, Theorem \ref{MAIN},  we show that $F_\varepsilon$ and $F_0$ subject to the constraint $u=0$ on $P_\delta$ share the same $\Gamma$-limit which is given by $$
\int_{\Omega} f_{\hom}(\nabla u)\, dx +\int_{\Omega} \varphi(u)\, dx$$ where $\varphi(z)$ is described by the capacitary formula induced by $F_{0}$ (see \eqref{capacita non lineare}) . It is worth noting that a similar analysis in the case of non local functionals of convolution type has been performed in a recent paper by R. Alicandro, M.S. Gelli and C. Leone  \cite{AGL}. From a technical point of view we mainly adopt the strategy exploited in \cite{AGL}, which in turn is inspired by \cite{AB}. The idea is to use a separation-of-scales argument, formalized in \cite{AB} and then deal with the non locality of our functionals to estimate the contribution near the perforations. A crucial role is played by a discrete version of the Gagliardo-Nirenberg-Sobolev inequality, which is derived from the corresponding non local variant proved in \cite{AGL}. This inequality allows us to show the convergence of minimum problems on unbounded domains, defining the approximating capacitary densities, to the limit energy density defined in \eqref{capacita non lineare}.\\
\vspace{0.3cm}

\noindent The paper is organized as follows. In Section $2$ we introduce some notation. In Section $3$ we present the setting of the problem and state the main result of the paper. In Section $4$ we recall some preliminary results. In Section $5$ we state and prove the discrete version of the Gagliardo-Nirenberg-Sobolev inequality and some other results which are instrumental  for the proof of the main theorem, which is the core of Section $6$.

\section{Notation}
In what follows $d,m \in \mathbb{N}$ will be two fixed natural numbers denoting the dimension of the reference and the target spaces of the functions we consider, respectively. The set of vectors $\{e_{1},\dots,e_{d}\}$ will denote the standard orthonormal basis in $\mathbb{R}^{d}$. Given $t \in \mathbb{R}$, $[t]$ denotes the integer part of $t$; for $\alpha \in \mathbb{Z}^{d}, r>0, Q(\alpha,r)=\alpha+(-r/2,r/2)^{d}$ (if $\alpha =0$, simply $Q_{r}$) is the open cube in $\mathbb{R}^{d}$ of center $\alpha$
and side length $r$. We denote by $\mathbb{S}^{d-1}$ the unit sphere in $\mathbb{R}^{d}$. If $A$ is a subset of $\mathbb{R}^{d}$ then $\dist(x,A)=\inf\{|y-x|:y \in A\}$; $\mathcal{A}^{\reg}(A)$ is the family of open subsets with Lipschitz boundary. We use standard notation for $\Gamma$-convergence \cite{dal}. Unless otherwise stated, $C$ will always denote a generic strictly positive constant that may change from line to line.

\section{Setting of the problem and the main result}
We fix $p \in (1,d)$ and we let $\Omega \subset \mathbb{R}^{d}$ be a bounded open set with Lipschitz boundary.
 For fixed $\varepsilon>0$, we denote by $\Omega_{\varepsilon}$ the lattice $\Omega_{\varepsilon}:=\varepsilon\mathbb{Z}^{d}\cap \Omega$ and later for a given infinitesimal sequence $\varepsilon_j$ we will use the notation $\Omega_{j}:=\Omega_{\varepsilon_j}$. We denote by $\mathcal{A}_{\varepsilon}(\Omega;\mathbb{R}^{m})$ the set of functions 
\begin{equation*}
\mathcal{A}_{\varepsilon}(\Omega;\mathbb{R}^{m}):=\{u:\Omega_{\varepsilon}\rightarrow \mathbb{R}^{m}\}.
\end{equation*}
We will identify the functions in $\mathcal{A}_{\varepsilon}(\Omega;\mathbb{R}^{m})$ by their piecewise constant interpolation on the cells of the lattice $\varepsilon\mathbb{Z}^{d}$ that is
\begin{equation*}
    \mathcal{A}_{\varepsilon}(\Omega;\mathbb{R}^{m})=\{u:\mathbb{R}^{d} \rightarrow \mathbb{R}^{m}: u\ \text{constant on}\  \alpha +[0,\varepsilon)^{d}\ \text{for any}\ \alpha \in \Omega_{\varepsilon}\}.
\end{equation*}
Given $\xi \in \mathbb{Z}^{d}$ and $E \subset \Omega$ we define 
\begin{equation}
    E_{\varepsilon}(\xi):=\{\alpha\in E |\ \alpha+\varepsilon\xi \in E\} \cap \varepsilon\mathbb{Z}^{d}.
    \label{set}
\end{equation}

\noindent Given a function $v \in \mathcal{A}_{\varepsilon}(\Omega;\mathbb{R}^{m})$, we denote by $D^{\xi}_{\varepsilon}v$ the different quotient along the direction $\xi$; i.e. for $\alpha \in \Omega_{\varepsilon}(\xi)$
\begin{equation}
    D^{\xi}_{\varepsilon}v(\alpha):=\frac{v(\alpha+\varepsilon \xi)-v(\alpha)}{\varepsilon|\xi|}
    \label{diffin}
\end{equation}
Let $\delta, r_{\delta}$ be given with $\delta>r_{\delta}>0$. We assume that 
\begin{equation*}
    \begin{split}
    &\delta=\delta_{\varepsilon}=N_{\varepsilon}\varepsilon \quad :N_{\varepsilon}\in \mathbb{N}\\
    &r_{\delta}=2n_{\varepsilon}\varepsilon \quad : n_{\varepsilon}\in \mathbb{N}.
\end{split}
\end{equation*}

\noindent With this setting we have that $\varepsilon\mathbb{Z}^{d}\supset \delta \mathbb{Z}^{d}$. For every $i\in \mathbb{Z}^{d}$, we define 
\begin{equation*}
    \mathcal{Q}(i\delta,r_{\delta}):=i\delta +\left(\left[-\frac{r_{\delta}}{2},\frac{r_{\delta}}{2}\right]^{d} \cap \varepsilon\mathbb{Z}^{d}\right)
\end{equation*}

\noindent and 

\begin{equation*}
    P_{\delta}:=\bigcup_{i\in \mathbb{Z}^{d}}\mathcal{Q}(i\delta,r_{\delta}).
\end{equation*}

\begin{figure}
    \centering
\includegraphics[width=0.5\linewidth]{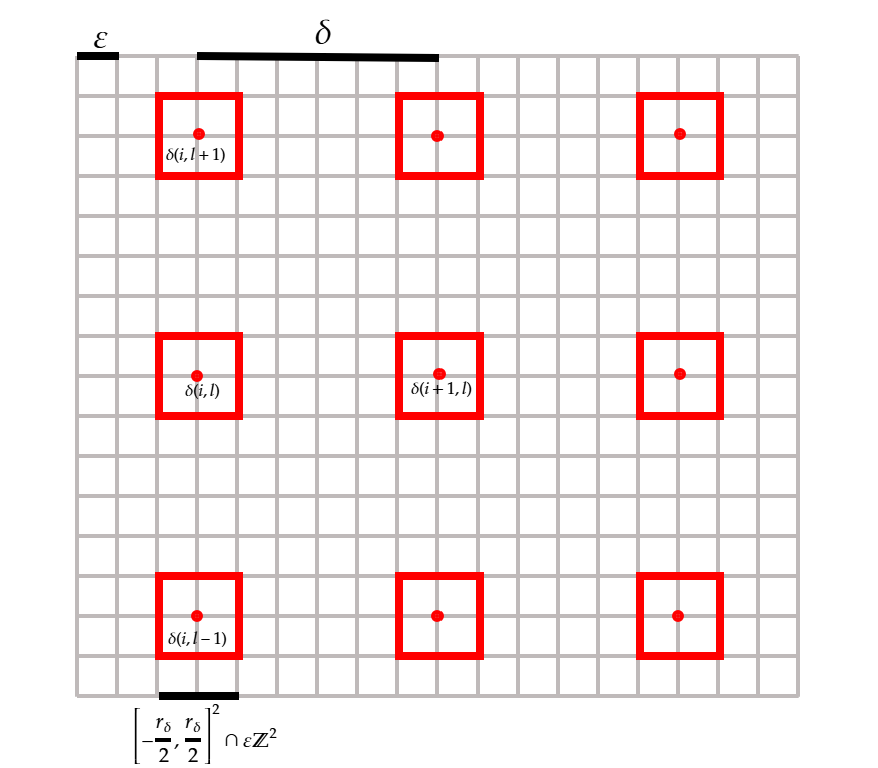}
    \caption{}
    \label{Figure1}
\end{figure}

\noindent Then, given $f:\mathbb{Z}^{d}\times \mathbb{R}^{m}\rightarrow[0,+\infty)$ we introduce the family of discrete functionals $F_{\varepsilon,\delta}:L^p(\Omega,\mathbb{R}^m)\to[0,+\infty]$ defined as follows

\begin{equation}
    F_{\varepsilon,\delta}(u):=
    \begin{cases}
    \displaystyle
        \sum_{\xi\in \mathbb{Z}^{d}}\sum_{\alpha \in \Omega_{\varepsilon}(\xi)}\varepsilon^{d}f(\xi,D^{\xi}_{\varepsilon}u(\alpha)) &\text{if}\ u\in\mathcal{A}_\varepsilon (\Omega;\mathbb{R}^m),\  u \equiv 0 \ \text{on} \ P_{\delta}\cap \Omega\\
        +\infty &\text{otherwise}.
    \end{cases}
    \label{1}
\end{equation}

\noindent We consider the following assumptions on $f$:
\begin{itemize}
    \item[(H)]$f(\xi,tz)=|t|^{p}f(\xi,z), \quad \forall t \in \mathbb{R},\ \forall (\xi,z) \in \mathbb{Z}^{d}\times \mathbb{R}^{m}$
    \item[(G)] For any $i=1,\dots, d$ and $\xi \in \mathbb{Z}^{d}$ the functions $m_{i}=\displaystyle\inf_{\zeta \in \mathbb{S}^{d-1}}f(e_{i},\zeta)$ and $M(\xi)=\displaystyle\sup_{\zeta \in \mathbb{S}^{d-1}}f(\xi,\zeta)$ satisfy
    \begin{itemize}
        \item[(G0)] there exists a constant $\lambda_{0}>0$ such that $\min\{m_{i}:i=1,\dots,d\}\geq \lambda_{0}$
        \item[(G1)] $\displaystyle\sum_{\xi\in \mathbb{Z}^{d}}M(\xi)<+\infty$ 
    \end{itemize}    \item[(L)]$|f(\xi,z_{1})-f(\xi,z_{2})|\leq CM(\xi)(|z_{1}|^{p-1}+|z_{2}|^{p-1})(|z_{1}-z_{2}|), \quad \forall z_{1}, z_{2} \in \mathbb{R}^{m}, \ \forall \xi \in \mathbb{Z}^{d}$
    
\end{itemize}
\begin{remark}
    From (G0) and (H) we get that for every $z \in \mathbb{R}^{m}$ and $i=1,\dots, d$, $f(e_{i},z)\geq \lambda_{0}|z|^{p}$. Indeed 
    \begin{equation*}
        f(e_{i},z)=|z|^{p}f\bigg(e_{i},\frac{z}{|z|}\bigg)\geq \lambda_{0}|z|^{p}.
    \end{equation*}
    Analogously for every $(\xi,z) \in \mathbb{Z}^{d}\times \mathbb{R}^{m}$, $f(\xi,z)\leq M(\xi)|z|^{p}$. Moreover, by (G1) we have that for every $\xi \in \mathbb{Z}^{d}$ and $\eta >0$ there exists $R_{\eta}>0$ such that
    \begin{equation*}
        \sum_{|\xi|>R_{\eta}}M(\xi)<\eta.
    \end{equation*}
    \label{oss1}
\end{remark}

\begin{remark}
    Hypothesis (L) trivially holds if  $f$ satisfies (H), (G), and $f(\xi,\cdot)$ is convex for every $\xi \in \mathbb{R}^{d}$.
\end{remark}
\noindent We also introduce the \textquoteleft truncated\textquoteright\  
 discrete functionals $F^T_{\varepsilon,\delta}:L^p(\Omega,\mathbb{R}^m)\to[0,+\infty]$ defined for every $T>0$ as 

\begin{equation}
    \label{funzionale troncato}F^{T}_{\varepsilon,\delta}(u):=
    \begin{cases}
    \displaystyle
        \sum_{\substack{\xi\in \mathbb{Z}^{d}\\
        |\xi|\leq T}}\sum_{\alpha \in \Omega_{\varepsilon}(\xi)}\varepsilon^{d}f(\xi,D^{\xi}_{\varepsilon}u(\alpha)) &\text{if}\ u\in\mathcal{A}_\varepsilon(\Omega;\mathbb{R}^m),\  u \equiv 0 \ \text{on} \ P_{\delta}\\
        +\infty &\text{otherwise}.
    \end{cases}
\end{equation}
  
\noindent Let us also consider the unconstrained family of functionals $\mathcal{F}_{\varepsilon}:L^{p}(\Omega;\mathbb{R}^{m})\rightarrow [0,+\infty]$ defined as
\begin{equation}
\mathcal{F}_{\varepsilon}(u)=
\begin{cases}
    \sum_{\xi \in \mathbb{Z}^{d}}\sum_{\alpha \in \Omega_{\varepsilon}(\xi)}\varepsilon^{d}f(\xi,D_{\varepsilon}^{\xi}u(\alpha)) &\text{if} \ u\in\mathcal{A}_\varepsilon(\Omega;\mathbb{R}^m),\\
    +\infty &\text{otherwise}
\end{cases}
\label{unconstrained functional}
\end{equation}

\noindent We introduce a localized version of such functionals: given an open set $A\subset \Omega$ we isolate the contribution due to interactions within $A$ by setting, for $u \in \mathcal{A}_{\varepsilon}(\Omega;\mathbb{R}^{m})$, 
\begin{equation}
    \mathcal{F}_{\varepsilon}(u,A):=\sum_{\xi \in \mathbb{Z}^{d}}\sum_{\alpha \in A_{\varepsilon}(\xi)}\varepsilon^{d}f(\xi,D_{\varepsilon}^{\xi}u(\alpha)).
    \label{10}
\end{equation}
The truncated version of \eqref{10} is defined, for any $T> 0$, as 
\begin{equation}
    \mathcal{F}_{\varepsilon}^{T}(u,A):=\sum_{\substack{\xi\in \mathbb{Z}^{d}\\
        |\xi|\leq T}}\sum_{\alpha \in A_{\varepsilon}(\xi)}\varepsilon^{d}f(\xi,D^{\xi}_{\varepsilon}u(\alpha)).
        \label{11}
\end{equation}
\noindent  
 In \cite[Theorem 4.1]{alicic}, it is proved the following $\Gamma$-convergence result.
\begin{theorem}
    Let $\mathcal{F}_{\varepsilon}$ be defined by (\ref{unconstrained functional}), with $f$ satisfying assumptions  (H) and (G). Then $(\mathcal{F}_{\varepsilon})$ $\Gamma$-converges with respect to the $L^{p}(\Omega;\mathbb{R}^{m})$-topology to the functional $\mathcal{F}:L^{p}(\Omega;\mathbb{R}^{m})\rightarrow [0,+\infty]$ defined as 
    \begin{equation}
    \label{Gamma limite AC}\mathcal{F}(u)=
        \begin{cases}
            \displaystyle
            \int_{\Omega}f_{\hom}(\nabla u)\, dx &\text{if}\ u \in W^{1,p}(\Omega;\mathbb{R}^{d}),\\
            +\infty &\textit{otherwise},
        \end{cases}
    \end{equation}
    where $f_{\hom}:\mathbb{R}^{m\times d}\rightarrow [0,+\infty)$ is given by the following homogenization formula:
    \begin{equation}
        f_{\hom}(M):=\lim_{h \rightarrow +\infty}\frac{1}{h^{d}}\min\bigg\{\sum_{\xi\in\mathbb{Z}^{d}}\sum_{\alpha \in (Q_{h})_{1}(\xi)}f(\xi,D_{1}^{\xi}u(\alpha)), \quad u \in \mathcal{A}_{1,M}(Q_{h})\bigg\}
        \label{13}
    \end{equation}
    with 
    \begin{equation*}
        \mathcal{A}_{1,M}(Q_{h}):=\{u \in \mathcal{A}_{1}(\mathbb{R}^{d};\mathbb{R}^{m}): u(\alpha)=M\alpha \  \text{if} \ (\alpha +[-1,1]^{d})\cap Q_{h}^{c}\neq \emptyset\}.
    \end{equation*}
    Moreover, there exist two positive constants $c<C$ such that, for every $M \in \mathbb{R}^{m \times d}$
    \begin{equation}
    \label{crescita dell' omogenizzata}
        c|M|^{p}\le f_{\hom}(M)\le C|M|^{p}.
    \end{equation}
    \label{51}
\end{theorem}

\noindent We now state the main result of the paper.

\begin{theorem}
    Let $F_{\varepsilon,\delta}$ be defined by \eqref{1}, with $f$ satisfying assumptions (H)-(G)-(L) and $1<p<d$. Assume moreover that 
    \begin{equation}
        \lim_{\delta \rightarrow 0}\frac{r_{\delta}}{\delta^{\frac{d}{d-p}}}=\gamma
    \end{equation}
    and that $\delta=\delta_{\varepsilon}$ is such that 
    \begin{equation}
        \lim_{\varepsilon\rightarrow 0}\frac{\varepsilon}{r_{\delta_{\varepsilon}}}=0
    \label{riscalamento taglia del reticolo e taglia della perforazione}\end{equation}
    for some real number $\gamma >0$. Then
    \begin{equation*}
        \Gamma(L^{p})-\lim_{\varepsilon\rightarrow 0}F_{\varepsilon,\delta_{\varepsilon}}(u)=
        \begin{cases}
        \displaystyle
            \int_{\Omega}f_{\hom}(\nabla u)\,dx+\gamma^{d-p}\int_{\Omega}\varphi(u)\, dx &if\ u \in W^{1,p}(\Omega;\mathbb{R}^{m}),\\
            +\infty &otherwise,
        \end{cases}
    \end{equation*}
    where $f_{hom}(M)$ is defined by (\ref{13}) and for every $z \in \mathbb{R}^{m}$
    \begin{equation}
    \label{capacita non lineare}
        \varphi(z):=\inf\bigg\{\int_{\mathbb{R}^{d}}f_{\hom}(\nabla u)\, dx: u \equiv 0\ \text{in} \ Q_{1}, u-z \in L^{p^{*}}(\mathbb{R}^{d};\mathbb{R}^{m}) \cap W^{1,p}_{\loc}(\mathbb{R}^{d};\mathbb{R}^{m}) \bigg\}.
    \end{equation}
    with $p^{*}=\frac{pd}{d-p}$.
    \label{MAIN}
\end{theorem}

\begin{remark}
    Taking into account the non degeneracy of the capacitary density \eqref{capacita non lineare} proved in Proposition \ref{non degeneracy of capacitary density} below, and arguing by comparison, one easily infers that the result stated in Theorem \ref{MAIN} can be "continuously" extended to the case $\gamma=0$ and $\gamma=+\infty$. More in details, if $\gamma=0$ and \eqref{riscalamento taglia del reticolo e taglia della perforazione} holds, then the functionals $F_{\varepsilon,\delta_{\varepsilon}}$ $\Gamma$-converge to the energy functional defined in \eqref{Gamma limite AC}, while, if $\gamma=+\infty$ and \eqref{riscalamento taglia del reticolo e taglia della perforazione} holds, the $\Gamma$-limit of the functionals $F_{\varepsilon,\delta_{\varepsilon}}$ is trivially $0$ if $u\equiv 0$ and $+\infty$ otherwise.
\end{remark}

\noindent For later use and reader's convenience we redefine the density functions we have introduced so far in the case $f(\xi,z)$ is replaced by $f^{T}(\xi,z)=\chi_{B_{T}}f(\xi,z)$. More precisely we set 
    \begin{equation*}
        {F}^{T}(u)=\int_{\Omega}f^{T}_{\hom}(\nabla u)\, dx+\gamma^{d-p}\int_{\Omega}\varphi^{T}(u)\, dx
    \end{equation*}
    with
    \begin{equation}
        f^{T}_{\hom}(M)=\lim_{h \rightarrow +\infty}\frac{1}{h^{d}}\min\bigg\{\sum_{\xi\in\mathbb{Z}^{d}}\sum_{\alpha \in (Q_{h})_{1}(\xi)}f^{T}(\xi,D_{1}^{\xi}u(\alpha)), \quad u \in \mathcal{A}_{1,M}(Q_{h})\bigg\}
        \label{f omogenizzata}
    \end{equation}
    and 
    \begin{equation}
        \varphi^{T}(z)=\inf\bigg\{\int_{\mathbb{R}^{d}}f^{T}_{\hom}(\nabla u)\, dx: u \equiv 0\ \text{in}\ Q_{1}, u-z \in L^{p^{*}}(\mathbb{R}^{d};\mathbb{R}^{m}) \cap W^{1,p}_{\loc}(\mathbb{R}^{d};\mathbb{R}^{m}) \bigg\}.
        \label{denista troncata}
    \end{equation}

\section{Preliminary results}
For every open set $A \subset \mathbb{R}^{d}$ and $u \in L^{p}(A;\mathbb{R}^{m})$ let us define
\begin{equation*}
G^{1,p}_{\varepsilon}(u,A):=\sum_{k=1}^{d}\sum_{\alpha\in A_{\varepsilon}(e_{k})}\varepsilon^{d}|D^{e_{k}}_{\varepsilon}u(\alpha)|^{p}
\end{equation*}
where $A_{\varepsilon}(e_{k})$ and $D^{e_{k}}_{\varepsilon}u(\alpha)$ are defined by (\ref{set}) and (\ref{diffin}), respectively. Concerning these type of energies, the following lemma proved in  \cite{alicic} is instrumental for our analysis.

\begin{lemma}[\cite{alicic}, Lemma 3.6]
    Let $A \in \mathcal{A}(\Omega)$, and set $\tilde{A}_{\varepsilon}:=\{x \in A: \dist(x,\partial A)>2\sqrt{d}\varepsilon\}.$ Then there exists a positive constant $C$ such that for any $\xi \in \mathbb{Z}^{d}$ and $u \in \mathcal{A}_{\varepsilon}(\Omega;\mathbb{R}^{m})$ there holds
    \begin{equation}
        \label{controllo long range with short range}
    \sum_{\alpha \in (\tilde{A}_{\varepsilon})_{\varepsilon}(\xi)}|D^{\xi}_{\varepsilon}u(\alpha)|^{p}\leq C\sum_{k=1}^{d}\sum_{\alpha \in A_{\varepsilon}(e_{k})}|D^{e_k}_{\varepsilon}u(\alpha)|^{p}.
    \end{equation}
    \label{AC}
\end{lemma}

\begin{remark}
\noindent As a direct consequence of the lemma above we can show that long range interactions are negligible for a sequence $u_{\varepsilon}$ such that $G_{\varepsilon}^{1,p}(u_{\varepsilon},\Omega)$ is uniformly bounded. Let us consider $\Omega'\subset \subset \Omega$ and $\eta>0$. Then,
\begin{equation*}
\begin{split}
    F_{\varepsilon}(u,\Omega')&=\sum_{\xi \in \mathbb{Z}^{d}}\sum_{\alpha\in (\Omega')_{\varepsilon}(\xi)}\varepsilon^{d}f(\xi,D^{\xi}_{\varepsilon}u(\alpha))\\
    &=\sum_{|\xi|\leq R_{\eta}}\sum_{\alpha \in (\Omega')_{\varepsilon}(\xi)}\varepsilon^{d}f(\xi,D^{\xi}_{\varepsilon}u(\alpha))+\sum_{|\xi|> R_{\eta}}\sum_{\alpha \in (\Omega')_{\varepsilon}(\xi)}\varepsilon^{d}f(\xi,D^{\xi}_{\varepsilon}u(\alpha)).
\end{split}
\end{equation*}
 Let us call $\tilde{\Omega}_{\varepsilon}:=\{x \in \Omega: \dist(x,\partial \Omega)>2\sqrt{d}\varepsilon\}$. It is clear that for $\varepsilon$ small $\Omega'\subset \subset \tilde{\Omega}_{\varepsilon}$. By (\ref{controllo long range with short range}) and Remark \ref{oss1} we get
 \begin{equation*}
     \begin{split}
         \sum_{|\xi|> R_{\eta}}\sum_{\alpha \in (\Omega')_{\varepsilon}(\xi)}\varepsilon^{d}f(\xi,D^{\xi}_{\varepsilon}u(\alpha))&\leq\sum_{|\xi|> R_{\eta}}\sum_{\alpha \in (\tilde{\Omega}_{\varepsilon})_{\varepsilon}(\xi)}\varepsilon^{d}f(\xi,D^{\xi}_{\varepsilon}u(\alpha))\\
     &\leq\sum_{|\xi|>R_{\eta}}M(\xi)\sum_{\alpha\in (\tilde{\Omega}_{\varepsilon})_{\varepsilon}(\xi)}\varepsilon^{d}|D^{\xi}_{\varepsilon}u(\alpha)|^{p}\\
     &\leq \sum_{|\xi|>R_{\eta}}M(\xi)\sum_{k=1}^{d}\sum_{\alpha \in  (\tilde{\Omega}_{\varepsilon})_{\varepsilon}(e_{k})}\varepsilon^{d}|D^{e_{k}}_{\varepsilon}u(\alpha)|^{p}\leq C\eta.
     \end{split}
 \end{equation*}
 Hence, by the arbitrariness of $\eta$, we get the thesis.
\end{remark}

\begin{remark}
    By Remark \ref{oss1} we deduce that, if $T>1$
    \begin{equation*}
    \mathcal{F}^{T}_{\varepsilon}(u,A)\geq \lambda_{0} G^{1,p}_{\varepsilon}(u,A).
    \end{equation*}
    Indeed,
    \begin{equation*}
    \begin{split}
        \mathcal{F}^{T}_{\varepsilon}(u,A)&=\sum_{|\xi|\leq T}\sum_{\alpha \in A_{\varepsilon}(\xi)}\varepsilon^{d}f(\xi,D^{\xi}_{\varepsilon}u(\alpha))\geq \sum_{k=1}^{d}\sum_{\alpha \in A_{\varepsilon}(e_{k})}\varepsilon^{d}f(e_{k},D^{e_{k
        }}_{\varepsilon}u(\alpha))\\
        &\geq \sum_{k=1}^{d}\sum_{\alpha \in A_{\varepsilon}(e_{k})}\varepsilon^{d}\lambda_{0}|D^{e_{k}}_{\varepsilon}u(\alpha)|^{p}=\lambda_{0}G^{1,p}_{\varepsilon}(u,A).
    \end{split}
    \end{equation*}
    \label{100}
\end{remark}

\noindent We recall the discrete version of the Poincarè-Wirtinger inequality and its rescaled version that have been proved in \cite  [Lemma 2, Lemma 3]{AS}.
\begin{lemma}
[\cite{AS}, Lemma 2]
Let $\Omega \subset \mathbb{R}^{d}$ a finite union of rectangles, and $p>1$. There exists $\varepsilon_{0}$ and a constant $C=C(p,\Omega)$ such that for all $\varepsilon<\varepsilon_{0}$ and $u:\Omega_{\varepsilon}\rightarrow \mathbb{R}^{m}$, having set $\Tilde{u}=\frac{1}{\sharp \Omega_{\varepsilon}}\sum_{\alpha \in \Omega_{\varepsilon}}u(\alpha)$, we have    \begin{equation}
        \sum_{\alpha \in \Omega_{\varepsilon}}\varepsilon^{d}|u(\alpha)-\tilde{u}|^{p}\leq CG^{1,p}_{\varepsilon}(u,\Omega).
        \label{PWd}
    \end{equation}
\label{53}
\end{lemma} 

\begin{lemma}[\cite{AS}, Lemma 3]
    \label{53bis}
    Let $\Omega, \varepsilon, p, C=C(p,\Omega)$ be as in Lemma \ref{53} and set $\Omega^{\delta}_{\varepsilon}=\delta \Omega\cap \varepsilon\mathbb{Z}^{d}$. Then, for $\varepsilon <\varepsilon_0$ and for all $u: \Omega^{\delta}_{\varepsilon}\rightarrow\mathbb{R}^{m}$ we have
    \begin{equation*}
        \sum_{\alpha \in \Omega^{\delta}_{\varepsilon}}\varepsilon^{d}|u(\alpha)-\tilde{u}|\leq C\delta^{p}\sum_{k=1}^{d}\sum_{\alpha \in \Omega_{\varepsilon}^{\delta}(e_{k})}\varepsilon^{d}\bigg|\frac{u(\alpha+e_{k}\varepsilon)-u(\alpha)}{\varepsilon}\bigg|^{p},
    \end{equation*}
    where $\tilde{u}=(\sharp\Omega^{\delta}_{\varepsilon})^{-1}\sum_{\alpha \in \Omega^{\delta}_{\varepsilon}} u (\alpha)$.
\end{lemma}
\begin{remark}
    The inequality \eqref{PWd} is still true if we consider the mean over a sufficiently regular set $E\subset \Omega$ with $|E|>0$,that is $\tilde{u}=\frac{1}{\sharp E_{\varepsilon}}\sum_{\alpha \in E_{\varepsilon}}u(\alpha)$. We note that in this case $\frac{\sharp \Omega_{\varepsilon}}{\sharp E_{\varepsilon}}\leq C$, where $C$ that does not depend on $\varepsilon$. \\
    Indeed
    \begin{equation*}
    \displaystyle\begin{split}
        \sum_{\alpha \in \Omega_{\varepsilon}}\bigg|u(\alpha)-\tilde{u}\bigg|^{p}&=\sum_{\alpha \in \Omega_{\varepsilon}}\bigg|u(\alpha)-\frac{1}{\sharp E_{\varepsilon}}\sum_{\beta \in E_{\varepsilon}}u(\beta)\bigg|^{p}
        =\sum_{\alpha \in \Omega_{\varepsilon}}\bigg|\frac{1}{\sharp E_{\varepsilon}}\sum_{\beta \in E_{\varepsilon}}u(\alpha)-u(\beta)\bigg|^{p}\\
        &\leq \sum_{\alpha \in \Omega_{\varepsilon}}\frac{1}{\sharp E_{\varepsilon}}\sum_{\beta \in E_{\varepsilon}}|u(\alpha)-u(\beta)|^{p}
        \leq \sum_{\alpha \in \Omega_{\varepsilon}}\frac{1}{\sharp E_{\varepsilon}}\sum_{\beta \in \Omega_{\varepsilon}}|u(\alpha)-u(\beta)|^{p}\\
        &=\sum_{\alpha \in \Omega_{\varepsilon}}\frac{\sharp \Omega_{\varepsilon}}{\sharp E_{\varepsilon}}\frac{1}{\sharp \Omega_{\varepsilon}}\sum_{\beta \in \Omega_{\varepsilon}}|u(\alpha)-u(\beta)|^{p}\leq C\sum_{\alpha \in \Omega_{\varepsilon}}\frac{1}{\sharp \Omega_{\varepsilon}}\sum_{\beta \in \Omega_{\varepsilon}}|u(\alpha)-u(\beta)|^{p}.
    \end{split}
    \end{equation*}
    Then we can proceed as in \cite{AS}.
   
    \label{OSSPW}
\end{remark}

\section{Supporting results}
In this section we state and prove some instrumental results for the proof of Theorem \ref{MAIN}.
\subsection{Gagliardo-Nirenberg-Sobolev type inequality}
In this subsection we state a crucial result for our analysis which may be of independent interest and resembles the Gagliardo-Nirenberg-Sobolev inequality in Sobolev spaces. Its proof is a direct consequence of the Gagliardo-Nirenberg-Sobolev type inequality for non local functionals \cite{AGL}, that we recall in the following.\\
\noindent For $r>0$
and given $p \geq 1$, let $T_{\varepsilon}:L^{p}(\mathbb{R}^{d};\mathbb{R}^{m})\rightarrow \mathcal{A}_{\tilde{r}\varepsilon}(\mathbb{R}^{d};\mathbb{R}^{m})$ be defined by
\begin{equation}
T_{\varepsilon}u(x):=\frac{1}{(\tilde{r}\varepsilon)^{d}}\int_{\tilde{r}\varepsilon k+[0,\tilde{r}\varepsilon)^{d}}u(y)\, dy \quad \ \text{on}\ \tilde{r}\varepsilon k+[0,\tilde{r}\varepsilon)^{d},\ k \in \mathbb{Z}^{d},
    \label{41}
\end{equation}
where $\tilde{r}:=\frac{r}{\sqrt{d+3}}$, and let us consider the functionals
    \begin{equation}
        \mathcal{G}^{r,p}_{\varepsilon}(u,A)=\int_{B_{r}}\int_{\mathbb{R}^{d}}|D^{\xi}_{\varepsilon}u(x)|^{p}\, dx\,d\xi.
        \label{Gtondo}
    \end{equation}
defined for every open set $A \subset \mathbb{R}^{d}$ and $u \in L^{p}(A;\mathbb{R}^{m})$

\begin{theorem}[\cite{AGL}, Theorem 5.1]
    Let $p \in [1,d)$. Then there exists a constant $C=C(p,d,r)>0$ such that for every $u \in L^{p}(\mathbb{R}^{d};\mathbb{R}^{m})$
    \begin{equation*}
    \bigg(\int_{\mathbb{R}^{d}}|T_{\varepsilon}u(x)|^{p^{*}}\, dx\bigg)^{\frac{p}{p^{*}}}\leq C\mathcal{G}^{r,p}_{\varepsilon}(u,\mathbb{R}^{d}),
    \end{equation*}
    where $p^{*}:=\frac{pd}{d-p}$.
    \label{Gagliardo}
\end{theorem}

\noindent Here we present a discrete version of the previous theorem. 

\begin{theorem}
    Let $p \in [1,d)$ and $p^{*}:=\frac{pd}{d-p}$. Then there exists a constant $C=C(p,d)>0$ such that for every $u \in \mathcal{A}_{\varepsilon}(\mathbb{R}^{d};\mathbb{R}^{m})\cap L^{p}(\mathbb{R}^{d},\mathbb{R}^{m})$

\begin{equation}
    \bigg(\int_{\mathbb{R}^{d}}|u(x)|^{p^{*}}\, dx\bigg)^{\frac{p}{p^{*}}}\leq CG^{1,p}_{\varepsilon}(u,\mathbb{R}^{d}).
        \label{44}
\end{equation}

\begin{proof}
We first observe that when $\tilde{r}=1$, which means $r=\sqrt{d+3}$, $T_{\varepsilon}u(x)=u(x)=u(\alpha)$ for $x \in \alpha+[0,\varepsilon)^{d} $. Thus, by Theorem \ref{Gagliardo} it is sufficient to prove that for $r=\sqrt{d+3}$
    \begin{equation*}
        \mathcal{G}^{r,p}_{\varepsilon}(u,\mathbb{R}^{d})\leq C G^{1,p}_{\varepsilon}(u,\mathbb{R}^{d}).
    \end{equation*}
    As consequence of \cite[Lemma 4.1]{AABPT}, applied with $\Omega=\mathbb{R}^{d}$, we have that for every $r>0$
    \begin{equation*}
        \mathcal{G}^{r,p}_{\varepsilon}(u,\mathbb{R}^{d})\leq C\mathcal{G}^{1,p}_{\varepsilon}(u,\mathbb{R}^{d}).
    \end{equation*}
    Hence, it is sufficient to prove that 
    \begin{equation*}
        \mathcal{G}^{1,p}_{\varepsilon}(u,\mathbb{R}^{d})\leq C G^{1,p}_{\varepsilon}(u,\mathbb{R}^{d}).
    \end{equation*}
    The energy $\mathcal{G}^{1,p}_{\varepsilon}(u,\mathbb{R}^{d})$ can be rewritten as 
\begin{equation*}
\mathcal{G}^{1,p}_{\varepsilon}(u,\mathbb{R}^{d})=\sum_{\alpha \in \varepsilon\mathbb{Z}^{d}}\int_{B_{1}}\int_{\alpha \in [0,\varepsilon)^{d}}\bigg|\frac{u(x+\varepsilon \xi)-u(x)}{\varepsilon}\bigg|^{p}\, dx\, d\xi.
\end{equation*}
Now, we fix $\xi \in B_{1}$ and $x \in \alpha +[0,\varepsilon)^{d}$. Then, either $x+\varepsilon\xi \in \alpha+[0,\varepsilon)^{d}$ and then $u(x+\varepsilon \xi)-u(x)=0$, or $x+\varepsilon\xi \in \alpha+\varepsilon e_{k}+[0,\varepsilon)^{d}$ for some $k=1,\dots,d$ and then $u(x+\varepsilon \xi)-u(x)=u(\alpha +\varepsilon e_{k})-u(\alpha)$. Now, let us define the set
\begin{equation}
\label{interfaccia}
    S^{k}_{\varepsilon,\alpha}:=\{x\in \alpha +[0,\varepsilon)^{d}: x+\varepsilon\xi \in \alpha +\varepsilon e_{k}+[0,\varepsilon)^{d}\},
\end{equation}
the interface in the $k$-th direction between two cells
\begin{equation*}
     J^{k}_{\varepsilon,\alpha}:=(\alpha+[0,\varepsilon]^{d})\cap (\alpha+\varepsilon e_{k}+[0,\varepsilon]^{d}),
\end{equation*}
and the stripe
\begin{equation*}
    J^{k,\xi}_{\varepsilon,\alpha}:=\{x+t\xi: -\varepsilon\leq t \leq \varepsilon,\ x \in J^{k}_{\varepsilon,\alpha}\}.
\end{equation*}

\noindent Since $S^{k}_{\varepsilon,\alpha}\subset J^{k,\xi}_{\varepsilon,\alpha}$, it follows that $|S^{k}_{\varepsilon,\alpha}|\leq C |J^{k,\xi}_{\varepsilon,\alpha}|\leq C\varepsilon^{d}|\xi|$. 
\noindent
Then we get  
\begin{equation*}
\begin{split}
\mathcal{G}^{1,p}_{\varepsilon}(u,\mathbb{R}^{d})&\leq C\int_{B_{1}}\sum_{k=1}^{d}\sum_{\alpha \in \varepsilon\mathbb{Z}^{d}}\int_{S^{k}_{\varepsilon,\alpha}}\bigg|\frac{u(x+\varepsilon \xi)-u(x)}{\varepsilon}\bigg|^{p}\, dx\, d\xi\\
&\leq
C\int_{B_{1}}|\xi|\, d\xi \sum_{k=1}^{d}\sum_{\alpha \in \varepsilon\mathbb{Z}^{d}}\varepsilon^{d}|D^{e_{k}}_{\varepsilon}u(\alpha)|^{p}\\
&=CG^{1,p}_{\varepsilon}(u,\mathbb{R}^{d}). 
\end{split}
\end{equation*}

\end{proof}
\label{Disuguaglianza di Gagliardo}
\end{theorem}

\subsection{Two technical lemmas}
In this section we recall two technical lemmas which will be used in the proof of Theorem \ref{MAIN} and have been proved in \cite{Siga}. The first one is the so called "joining lemma", which allows to reduce our $\Gamma$-convergence to sequences of converging functions that are constants on suitable annuli surrounding the perforations. The second one is the "truncation lemma" that, by the composition with a suitable Lipschitz function, allows us to replace a given sequence $u_{j}$ with equibounded energies and $L^{p}$-norms, by a new sequence uniformly bounded in $L^{\infty}$, with a small gap in energy. 
\begin{lemma}[\cite{Siga}, Lemma 6.1]

Let $\{u_j\}_j$ be a sequence such that $u_j \in \mathcal{A}_{\varepsilon_j}(\Omega;\mathbb{R}^{m})$ and $u_j\rightarrow u$ in $L^{p}(\Omega;\mathbb{R}^{m})$ for some $u \in W^{1,p}(\Omega;\mathbb{R}^{m})$. Assume that 
\begin{equation*}
\sup_{j}\sum_{\substack{\xi \in \mathbb{Z}^{d}\\ |\xi|\leq T}}\sum_{\alpha \in \Omega_{\varepsilon}(\xi)}\varepsilon_{j}^{d}f(\xi,D^{\xi}_{\varepsilon_{j}}u_{j}(\alpha))<+\infty.
\end{equation*}

\noindent Let $(\rho_{j})$ be a sequence of the form $\rho_{j}=\beta\delta_{j}$ where $\beta <\frac{1}{2}$ with $\varepsilon_{j}=o(\rho_{j})$ as $j\rightarrow +\infty$.
 We denote by $Z_{j}(\Omega)$ the set of indices 
\begin{equation}
    Z_{j}(\Omega)=\left\{i \in \mathbb{Z}^{d}\cap \frac{\Omega}{\delta_j}: \dist(i\delta_{j},\partial\Omega)>\delta_{j}\right\}.
    \label{insieme di indici}
\end{equation}
Let $k \in \mathbb{N}$ be fixed. Then for all $i \in Z_{j}(\Omega)$ there exists $k_{i} \in \{0,\dots, k-1\}$ such that, having set  
\begin{equation*}
C^{i}_{j}:=Q\bigg(i\delta_{j},\bigg[\frac{\rho_{j}}{\varepsilon_{j}}2^{-k_{i}}\bigg]\varepsilon_{j}\bigg)\bigg\backslash  Q\bigg(i\delta_{j},\bigg[\frac{\rho_{j}}{\varepsilon_{j}}2^{-k_{i}-1}\bigg]\varepsilon_{j}\bigg)
\end{equation*}

\begin{equation*}
    \rho^{i}_{j}:=\bigg[\frac{3}{4}\frac{\rho_{j}}{\varepsilon_{j}}2^{-k_{i}}\bigg]\varepsilon_{j},
\end{equation*}

\begin{equation*}
     u_{j}^{i}:=\frac{1}{\sharp C^{i}_{j}}\sum_{\alpha \in C^{i}_{j}} u_{j}(\alpha),
\end{equation*}
there exists a sequence $w_{j} \in \mathcal{A}_{\varepsilon_{j}}(\Omega;\mathbb{R}^{m})$ such that $w_{j} \rightarrow u$ in $L^{p}(\Omega;\mathbb{R}^{m})$ and 

\begin{equation}
    w_{j}=u_{j}\  \text{on}\ \Omega_{j} \setminus \bigcup_{i \in Z_{j}(\Omega)}C_{j}^{i},
    \label{outcornici}
\end{equation}

\begin{equation}
    w_{j}=u_{j}^{i}\ \text{on}\ Q(i\delta_j,\rho^{i}_{j}+\varepsilon T)\setminus Q(i\delta_j,\rho^{i}_{j}-\varepsilon T),
    \label{bordointerno}
\end{equation}

\begin{equation}
    \bigg|\sum_{|\xi|\leq T}\sum_{\alpha \in \Omega_{\varepsilon_j}(\xi)}\varepsilon_{j}^{d}(f(\xi,D^{\xi}_{\varepsilon_{j}}u_{j}(\alpha))-f(\xi,D^{\xi}_{\varepsilon_{j}}w_{j}(\alpha)))\bigg|\leq \frac{C}{k}.
    \label{diffenergia}
\end{equation}
   
\label{JL}
\end{lemma}

\begin{remark}
In \cite[Lemma 6.1]{Siga} condition \eqref{bordointerno} is replaced by $w_{j}=u_{j}^{i}$ on $\partial Q(i\delta_j,\rho^{i}_{j})$, however it can be easly inferred from that proof that $w_j$ can be chosen so that \eqref{bordointerno} is satisfied. 
\end{remark}

\begin{lemma}[\cite{Siga}, Lemma 7.2]
    Let $\{u_j\}_j$ be a sequence such that $u_{j} \in \mathcal{A}_{\varepsilon_j}(\Omega;\mathbb{R}^{m})$ with $\sup_{j}(\mathcal{F}_{\varepsilon}(u_{j})+\|u_{j}\|_{L^{p}(\Omega;\mathbb{R}^{m})})<+\infty$. Then for every $\eta >0$ and $L \in \mathbb{N}$ there exist a subsequence (not relabeled), a constant $R_{L}>L$ and a Lipschitz function $t_{L}:\mathbb{R}^{m}\rightarrow \mathbb{R}^{m}$ with Lipschitz constant $1$ such that $t_{L}(z)=z$ if $|z|<R_L$,  $t_{L}(z)=0$ if $|z|>2R_{L}$ and it holds
    \begin{equation*}
        \liminf_{j}\mathcal{F}_{\varepsilon_{j}}(t_{L}(u_{j}))\leq \liminf_{j}\mathcal{F}_{\varepsilon_{j}}(u_{j})+\eta.
    \end{equation*}
    \label{TL}
\end{lemma}

\noindent We recall a compactness result in the strong $L^{p}$ topology for sequences of functions whose energy defined in \eqref{Gtondo} is uniformly bounded.

\begin{theorem}[\cite{AABPT}, Theorem 4.2]
\label{compactness for concolution}
    Let $A$ be any open set of $\mathbb{R}^{d}$ with bounded Lipschitz boundary. Let $\{u_{\varepsilon}\}_{\varepsilon}\subset L^{p}(A;\mathbb{R}^{m})$ be such that for some $r>0$
    \begin{equation*}
        \sup_{\varepsilon>0}\{\|u_{\varepsilon}\|_{L^{p}(A;\mathbb{R}^{m})}+\mathcal{G}^{r,p}_{\varepsilon}(u_{\varepsilon},A)\}<+\infty.
    \end{equation*}
    Then, given $\varepsilon_j\rightarrow 0, \{u_{\varepsilon_j}\}_{j}$ is relatively compact in $L^{p}(A;\mathbb{R}^{m})$ and every limit of a converging subsequence is in $W^{1,p}(A;\mathbb{R}^{m})$.
\end{theorem}
\subsection{Approximating capacitary energy densities}

In this subsection we introduce and investigate the main properties of suitable energy densities
defined through minimum problems of capacitary type.\\
\noindent
For any $\varepsilon >0$, $T> 2$, $R>1 + \varepsilon T$, $z \in \mathbb{R}^{m}$, set
\begin{equation}
\varphi_{\varepsilon,T,R}(z):=\inf\{\mathcal{F}^{T}_{\varepsilon}(v,Q_{R}): v \in \mathcal{A}_{\varepsilon,T,z}(Q_{R};\mathbb{R}^{m})\}
    \label{20}
\end{equation}
where $\mathcal{F}^{T}_\varepsilon$ is defined as in \eqref{11} and 
\begin{equation*}
\mathcal{A}_{\varepsilon,T,z}(Q_{R};\mathbb{R}^{m}):=\{v\in \mathcal{A}_{\varepsilon}(Q_{R};\mathbb{R}^{m}): v \equiv 0\ \text{in} \ Q_{1}, v \equiv z \ \text{on} \ \partial^{\varepsilon T}Q_{R}\},
\end{equation*}
with the notation $\partial^{\varepsilon T}Q_{R}=Q_{R+\varepsilon T}\setminus Q_{R-\varepsilon T }$.

\noindent
 Given $v \in \mathcal{A}_{\varepsilon,T,z}(Q_{R};\mathbb{R}^{m})$, consider the function $u:=v-z$. Hence $u\equiv -z$ in $Q_{1}$ and $u \equiv 0$ on $\partial^{\varepsilon T}Q_{R}$. Then we can extend $u$ to the whole $\mathbb{R}^{d}$ by setting
\begin{equation*}
    u \equiv 0 \ \text{on} \ \mathbb{R}^{d}\setminus Q_{R-\varepsilon T}.
\end{equation*}
Identifying  $u$ with its extension, it is clear that $\mathcal{F}^{T}_{\varepsilon}(u,Q_{R})=\mathcal{F}^{T}_{\varepsilon}(u,\mathbb{R}^{d})$.
Hence (\ref{20}) can be written as
\begin{equation}
\varphi_{\varepsilon,T,R}(z)=\inf\{\mathcal{F}^{T}_{\varepsilon}(u, \mathbb{R}^{d}): u \in  \mathcal{A}_{\varepsilon}(\mathbb{R}^{d};\mathbb{R}^{m}), u \equiv -z \ \text{in} \ Q_{1}, u \equiv 0 \ \text{in} \ \mathbb{R}^{d}\setminus Q_{R-\varepsilon T}\}.
    \label{22}
\end{equation}

\begin{remark}
    If $f$ satisfies (G), (H) and (L), by applying the direct methods of Calculus of Variations, the infimum problem defining  $\varphi_{\varepsilon,T,R}(z)$ in \eqref{20} is in fact a minimum.
\end{remark}

\noindent Next result establishes growth conditions of order $p$ for $\varphi_{\varepsilon, T, R}(\cdot)$. 

\begin{proposition}
    Let $f$ satisfy assumptions $(H)$ and $(G)$, and let $T>2$, $\varepsilon_{0}>0$, and $R$ be fixed such that $R-2\varepsilon_{0}T\geq 2$. Then, for every $0<\varepsilon\leq \varepsilon_{0}$ there exists $c_{1}, c_{2}>0$ such that
    \begin{equation}
        c_{1}|z|^{p}\leq \varphi_{\varepsilon, T, R}(z) \quad \forall z \in \mathbb{R}^{m},
        \label{25}
    \end{equation}
    \begin{equation}
    \varphi_{\varepsilon, T, R}(z)\leq c_{2}|z|^{p} \quad \forall z \in \mathbb{R}^{m}.
        \label{26}
    \end{equation}
    
    \begin{proof}
    We firstly prove inequality (\ref{25}). The proof relies on a suitable lower bound of $\mathcal{F}^{T}_{\varepsilon}(u,Q_{R})$. In order to avoid many technicalities we detail the proof in the case $d=2$. We will associate to any admissible function $u$ a piecewise affine function obtained by linearly interpolating the values of $u$ on a triangulation of the lattice $\mathbb{Z}^{2}$. The argument can be generalized to any dimension by linearly interpolating the values of $u$ on a partition of the cells of the lattice $\mathbb{Z}^{d}$ in $d$-simplices, known in literature as the Kuhn decomposition (see for example \cite{ALP}). \\ Set
    \begin{equation*}
        T^{-}:=\{(x_{1},x_{2})\in [0,1]^{2}: x_{2}\leq 1-x_{1}\},
    \end{equation*}
    
    \begin{equation*}
         T^{+}:=\{(x_{1},x_{2})\in [0,1]^{2}: x_{2}\geq 1-x_{1}\},
    \end{equation*}
    and, given any admissible function $u$ in (\ref{22}),
    let us define the piecewise-affine interpolation of
    $u$ on the cells of the lattice as follows
    \begin{equation*}
    \hat{u}(x):=u(\alpha)+\frac{u(\alpha+\varepsilon e_{1})-u(\alpha)}{\varepsilon}(x_{1}-\alpha_{1})+\frac{u(\alpha+\varepsilon e_{2})-u(\alpha)}{\varepsilon}(x_{2}-\alpha_{2})\ \quad \text{if}\ x \in \alpha +\varepsilon T^{-}
    \end{equation*}

    \begin{equation*}
    \begin{split}
        \hat{u}(x)=u(\alpha+\varepsilon(e_1+e_2))&+\frac{u(\alpha+\varepsilon(e_1+e_2))-u(\alpha+\varepsilon e_2)}{\varepsilon}(x_1-\alpha_1-\varepsilon)\\
        &+\frac{u(\alpha+\varepsilon(e_1+e_2))-u(\alpha+\varepsilon e_1)}{\varepsilon}(x_2-\alpha_2-\varepsilon)\ \quad \text{if}\  x\in \alpha +\varepsilon T^{+}
    \end{split}
    \end{equation*}
    Since $u\equiv -z$ in $Q_{1}$ and $u\equiv 0$ in $\mathbb{R}^{d}\setminus Q_{R-\varepsilon T}$ we infer that, for $\varepsilon$ small enough $\hat{u} \equiv -z$ in $Q_{\frac{1}{2}}$ and $\hat{u} \equiv 0$ on $\mathbb{R}^{d}\setminus Q_{R}$. Let us call $P^{\varepsilon,1}_{\alpha}:=\alpha+\varepsilon P_1$ and $P^{\varepsilon,2}_{\alpha}:=\alpha+\varepsilon P_2$ where 
    \begin{equation*}
        P_1=\{(x_1,x_2) \in \mathbb{R}^{2}: 0\leq x_{1}\leq 1, x_{1}\leq x_{2}\leq 1-x_{1}\},
    \end{equation*}
    \begin{equation*}
        P_2=\{(x_1,x_2) \in \mathbb{R}^{2}: 0\leq x_{2}\leq 1, -x_{2}\leq x_{1}\leq 1-x_{2}\}.
    \end{equation*}
    One can easily show that
    \begin{equation*}
        \partial_{x_{k}}\hat{u}(x)=D^{e_{k}}_{\varepsilon}u(\alpha), \quad k=1,2\ \text{and}\ x\in  P^{\varepsilon,k}_{\alpha}.
    \end{equation*}
    
     \noindent Since $|P^{\varepsilon,1}_{\alpha}|=|P^{\varepsilon,2}_{\alpha}|=\varepsilon^{2}$, we get the following estimate
    \begin{equation*}
        \begin{split}
            \int_{\mathbb{R}^{2}}|\nabla \hat{u}|^{p}&\leq C\bigg(\int_{\mathbb{R}^{2}}|\partial_{x_{1}}\hat{u}|^{p}+|\partial_{x_{2}}\hat{u}|^{p}\bigg)\\
            &=C\sum_{\alpha \in \varepsilon \mathbb{Z}^{2}}\bigg(\int_{P^{\varepsilon,1}_{\alpha}}|\partial_{x_{1}}\hat{u}|^{p}+\int_{P^{\varepsilon,2}_{\alpha}}|\partial_{x_{2}}\hat{u}|^{p}\bigg)\\
            &=C\sum_{\alpha \in \varepsilon \mathbb{Z}^{2}}\varepsilon^{2}(|D^{e_{1}}_{\varepsilon}u(\alpha)|^{p}+|D^{e_{2}}_{\varepsilon}u(\alpha)|^{p})\\
            &=C\sum_{k=1}^{2}\sum_{\alpha \in \varepsilon \mathbb{Z}^{2}}\varepsilon^{2}|D^{e_{k}}_{\varepsilon}u(\alpha)|^{p}=CG^{1,p}_{\varepsilon}(u,\mathbb{R}^{2}).
        \end{split}
    \end{equation*}
     By Remark \ref{100} if follows
     \begin{equation*}
    \mathcal{F}^{T}_{\varepsilon}(u,\mathbb{R}^{2})\geq \lambda_{0}G_{\varepsilon}^{1,p}(u,\mathbb{R}^{2})
         \geq C\int_{\mathbb{R}^{2}}|\nabla \hat{u}|^{p}
         \geq |z|^{p}\capi(Q_{\frac{1}{2}},\mathbb{R}^{2})
         =c_{1}|z|^{p}.
     \end{equation*}

    \noindent
        We now prove inequality (\ref{26}). Consider a function $u \in C^{1}_{c}(Q_{R-2\varepsilon T};\mathbb{R}^{m})$ such that $u\equiv -z$ in $Q_{\frac{3}{2}}$, and for every $\alpha \in \varepsilon\mathbb{Z}^{d}$ let us define the discrete function
        \begin{equation*}
            \hat{u}(\alpha)=\frac{1}{\varepsilon^{d}}\int_{\alpha +[0,\varepsilon)^{d}}u(x)\, dx.
        \end{equation*}
         Notice that $\hat{u}$ is admissible for the problem defining $\varphi_{\varepsilon, T, R}$ in \eqref{22}. We show that there exists a constant $C$ such that, for $\varepsilon \leq \varepsilon_{0}$
        
        \begin{equation}
        G_{\varepsilon}^{1,p}(\hat{u},\mathbb{R}^{d})\leq C \int_{\mathbb{R}^{d}}|\nabla u|^{p}\, dx.
            \label{27}
        \end{equation}

        \noindent Indeed, fixed $k=1,\dots, d$ and $x \in \mathbb{R}^{d}$, we have
        
        \begin{equation}
        \frac{u(x+\varepsilon e_{k})-u(x)}{\varepsilon}= \int_{0}^{1}\partial_{x_{k}}u(x+\varepsilon e_{k}t)\, dt.
           \label{stima del differenze finite con il gradiente}
        \end{equation}
        By Jensen's inequality, Fubini's Theorem and (\ref{stima del differenze finite con il gradiente}) we have the following estimate
        \begin{equation*}
            \begin{split}
                G^{1,p}_{\varepsilon}(\hat{u},\mathbb{R}^{d})&=\sum_{k=1}^{d}\sum_{\alpha \in \varepsilon\mathbb{Z}^{d}}\varepsilon^{d}|D^{e_k}_{\varepsilon}\hat{u}(\alpha)|^{p}\\
                &=\sum_{k=1}^{d}\sum_{\alpha \in \varepsilon\mathbb{Z}^{d}}\varepsilon^{d}\left|\frac{1}{\varepsilon^{d}}\left[\int_{\alpha +[0,\varepsilon)^{d}}\frac{u(x+\varepsilon e_k)-u(x)}{\varepsilon}\, dx\right]\right|^{p}\\
                &\leq \sum_{k=1}^{d}\sum_{\alpha \in \varepsilon\mathbb{Z}^{d}} \int_{\alpha +[0,\varepsilon)^{d}}\bigg|\frac{u(x+\varepsilon e_k)-u(x)}{\varepsilon}\bigg|^{p}\, dx\\
                &\leq \sum_{k=1}^{d}\sum_{\alpha \in \varepsilon\mathbb{Z}^{d}} \int_{\alpha +[0,\varepsilon)^{d}} \int_{0}^{1}|\nabla u(x+\varepsilon e_k t)|^{p}  \, dt\, dx \\
                &=\sum_{k=1}^{d}\sum_{\alpha \in\varepsilon\mathbb{Z}^{d}}\int_{0}^{1}\int_{\alpha +[0,\varepsilon)^{d}} |\nabla u(x)|^{p}\, dx\, dt\\
                &\leq C\int_{\mathbb{R}^{d}}|\nabla u(x)|^{p}\, dx.
                \end{split}
        \end{equation*}
   
        \noindent By (G1), Lemma \ref{AC} and \eqref{27}  we get
        \begin{equation*}
            \begin{split}
                 \mathcal{F}^{T}_{\varepsilon}(\hat{u},\mathbb{R}^{d})
                &\leq \sum_{\substack{\xi \in \mathbb{Z}^{d}\\|\xi|\leq T}}\sum_{\alpha \in \varepsilon\mathbb{Z}^{d}}\varepsilon^{d}M(\xi)|D^{\xi}_{\varepsilon}\hat{u}(\alpha)|^{p}\\
                &\leq C \sum_{\xi\in \mathbb{Z}^{d}}M(\xi) G^{1,p}_{\varepsilon}(\hat{u},\mathbb{R}^{d})\\
                &
                \leq C \sum_{\xi \in \mathbb{Z}^{d}}M(\xi)\int_{\mathbb{R}^{d}}|\nabla u(x)|^{p}\, dx.
                \end{split}
                \end{equation*}

                    \noindent Hence, by \eqref{22} and a density argument 

                    \begin{equation*}
            \begin{split}
                \varphi_{\varepsilon, T, R}(z)&\leq C \sum_{\xi\in \mathbb{Z}^{d}}M(\xi) \inf \bigg\{ \int_{Q_{R-2\varepsilon T}}|\nabla u|^{p}\, dx: u \equiv -z \ \text{in}\ Q_{1}, u \in C^{1}_{c}(Q_{R-2\varepsilon T};\mathbb{R}^{m})\bigg\}\\
                &=C \sum_{\xi\in \mathbb{Z}^{d}}M(\xi) \capi\displaystyle_{p}(Q_{1},Q_{R-2\varepsilon T})|z|^{p}\leq C \sum_{\xi\in \mathbb{Z}^{d}}M(\xi) \capi\displaystyle_{p}(Q_{1},Q_{2})|z|^{p}=c_{2}|z|^{p}.
            \end{split}
        \end{equation*}
    \end{proof}
    \label{non degeneracy of capacitary density}
\end{proposition}

\noindent In the next proposition, we show that the functions $\varphi_{\varepsilon, T, R}$ are uniformly locally Lipschitz continuous. 

\begin{proposition}
\label{200}
    Let $f$ satisfy assumptions (H), (G), (L). Then there exists a constant $C>0$ independent of $\varepsilon, T$ and $R$ such that for every $z, w \in \mathbb{R}^{m}$ we have
    \begin{equation}
        \label{30}
        |\varphi_{\varepsilon,T,R}(w)-\varphi_{\varepsilon,T,R}(z)|\leq C(|z|^{p-1}+|w|^{p-1})|w-z|.
    \end{equation}
      \begin{proof}
        Note that by (\ref{26}), if $z=0$ or $w=0$ the inequality is trivial. Then, we can suppose $z$ and $w$ both not null and consider the map $\phi: \mathbb{R}^{m}\rightarrow \mathbb{R}^{m}$ defined by
        \begin{equation*}
            \phi(\zeta)=\frac{|w|}{|z|}\mathcal{R}^{w}_{z}(\zeta),
        \end{equation*}
        where $\mathcal{R}^{w}_{z}$ is a rotation that maps $z$ into $\frac{|z|}{|w|}w$. Note that $\phi(0)=0$, $\phi(z)=w$ and 
        \begin{equation}
            \|\nabla \phi\|_{\infty}\leq C\frac{|w|}{|z|}, \quad \|\nabla \phi -I\|_{\infty}\leq C \frac{|w-z|}{|z|}.
            \label{31}
        \end{equation}
        \noindent
        Let $v_{z} \in \mathcal{A}_{\varepsilon,T,z}(Q_{R};\mathbb{R}^{m})$ be such that 
        \begin{equation}
        \mathcal{F}^{T}_{\varepsilon}(v_{z},Q_{R})= \varphi_{\varepsilon,T,R}(z).
            \label{32}
        \end{equation}
        Now set 
        \begin{equation*}
            v_{w}:=\phi \circ v_{z},
        \end{equation*}
        and observe that $v_{w} \in \mathcal{A}_{\varepsilon,T,w}(Q_{R};\mathbb{R}^{m})$. Hence, by (\ref{32})
        \begin{equation*}
                \varphi_{\varepsilon,T,R}(w)\leq \mathcal{F}^{T}_{\varepsilon}(v_{w},Q_{R})=\varphi_{\varepsilon,T,R}(z)+\mathcal{F}^{T}_{\varepsilon}(v_{w},Q_{R})-\mathcal{F}^{T}_{\varepsilon}(v_{z},Q_{R}).
        \end{equation*}
        By hypothesis (L) and (\ref{31})
        \begin{equation*}
            \begin{split}
                |f(\xi,D^{\xi}_{\varepsilon}v_{w}(\alpha))&-f(\xi,D^{\xi}_{\varepsilon}v_{z}(\alpha))|\leq CM(\xi)(|D^{\xi}_{\varepsilon}v_{z}(\alpha)|^{p-1}+|D^{\xi}_{\varepsilon}v_{w}(\alpha)|^{p-1})|D^{\xi}_{\varepsilon}v_{z}(\alpha)-D^{\xi}_{\varepsilon}v_{w}(\alpha)|\\
                &=CM(\xi)(|D^{\xi}_{\varepsilon}v_{z}(\alpha)|^{p-1}+|\nabla \phi \cdot D^{\xi}_{\varepsilon}v_{z}(\alpha)|^{p-1})|D^{\xi}_{\varepsilon}v_{z}(\alpha)[\nabla \phi-I]|\\
                &\leq CM(\xi)(|D^{\xi}_{\varepsilon}v_{z}(\alpha)|^{p-1}+\|\nabla \phi\|_{\infty}^{p-1}|D^{\xi}_{\varepsilon}v_{z}(\alpha)|^{p-1})\|\nabla \phi-I\|_{\infty}|D^{\xi}_{\varepsilon}v_{z}(\alpha)|\\
                &\leq CM(\xi)|D^{\xi}_{\varepsilon}v_{z}(\alpha)|^{p-1}\left(1+C\frac{|w|^{p-1}}{|z|^{p-1}}\right)\frac{|w-z|}{|z|}|D^{\xi}_{\varepsilon}v_{z}(\alpha)|\\
                &=CM(\xi)\frac{|z|^{p-1}+|w|^{p-1}}{|z|^{p}}|w-z||D^{\xi}_{\varepsilon}v_{z}(\alpha)|^{p}
            \end{split}
        \end{equation*}
        Then, 
        \begin{equation}
            \varphi_{\varepsilon,T,R}(w)\leq \varphi_{\varepsilon,T,R}(z)+C\frac{|z|^{p-1}+|w|^{p-1}}{|z|^{p}}|w-z|\sum_{|\xi|\leq T}M(\xi)\sum_{\alpha \in (Q_{R})_{\varepsilon}(\xi)}|D^{\xi}_{\varepsilon}v_{z}(x)|^{p}.
            \label{33}
        \end{equation}
        By (G0), (\ref{32}) and (\ref{26})
        \begin{equation}
            G^{1,p}_{\varepsilon}(v_{z},Q_{R})\leq C\mathcal{F}^{T}_{\varepsilon}(v_{z},Q_{R})=C\varphi_{\varepsilon,T,R}(z)\leq C|z|^{p}.
        \label{bound di G}\end{equation}
        Observe that $G_{\varepsilon}^{1,p}(v_{z},Q_{R})=G_{\varepsilon}^{1,p}(v_{z},\mathbb{R}^{d})$. Then, by Lemma \ref{AC}
        \begin{equation*}
           \sum_{\alpha \in (Q_{R})_{\varepsilon}(\xi)}|D^{\xi}_{\varepsilon}v_{z}(\alpha)|^{p}\leq CG^{1,p}_{\varepsilon}(v_{z},Q_{R}).
        \end{equation*}
        Hence,
        by (\ref{33}) and (\ref{bound di G}) we get
        \begin{equation*}
            \varphi_{\varepsilon,T,R}(w)\leq \varphi_{\varepsilon,T,R}(z)+C(|z|^{p-1}+|w|^{p-1})|w-z|.
        \end{equation*}
        Reversing the role of $z$ and $w$, we get (\ref{30}).
      \end{proof}
\end{proposition}

\subsection{Asymptotics of the approximating capacitary energy density}

We now show that, if $R_{\varepsilon}\rightarrow +\infty$ as $\varepsilon \rightarrow 0$ the functions $\varphi_{\varepsilon,T,R_{\varepsilon}}(z)$ approximate the energy density $\varphi^{T}(z)$. Remember that they are defined respectively as
\begin{equation}
    \label{50}
    \varphi^{T}(z):=\inf\bigg\{\int_{\mathbb{R}^{d}}f^{T}_{\hom}(\nabla v)\, dx : v-z \in L^{p^{*}}(\mathbb{R}^{d};\mathbb{R}^{m}), v \equiv 0 \ \text{in} \ Q_{1}, v \in W^{1,p}_{\loc}(\mathbb{R}^{d};\mathbb{R}^{m})\bigg\},
\end{equation}
\begin{equation}
    \label{201}
    \varphi_{\varepsilon,T,R_{\varepsilon}}(z):=\inf\big\{\mathcal{F}^{T}_{\varepsilon}(u,\mathbb{R}^{d}): u \in \mathcal{A}_{\varepsilon}(\mathbb{R}^{d};\mathbb{R}^{m}), u \equiv -z\ \text{in}\ Q_{1}, u \equiv 0\ \text{in}\ \mathbb{R}^{d}\setminus Q_{R_{\varepsilon}-\varepsilon T}\big\}.
\end{equation}

\begin{proposition}
    Let $f$ satisfy assumptions (H), (G) and (L).
    Let $\varphi^{T}$ and $\varphi_{\varepsilon,T,R_{\varepsilon}}$ be defined by (\ref{50}) and (\ref{201}). Then, if $R_{\varepsilon}\rightarrow +\infty$ as $\varepsilon \rightarrow 0$ it holds
    \begin{equation}
        \label{202}
        \lim_{\varepsilon\rightarrow 0}\varphi_{\varepsilon,T,R_{\varepsilon}}(z)=\varphi^{T}(z)
    \end{equation}
    uniformly on compact sets.
    \begin{proof}
        By Proposition \ref{200} it is sufficient to prove that (\ref{202}) holds pointwise. We will show that \eqref{202}  is a consequence of Theorem \ref{51} and Theorem \ref{Disuguaglianza di Gagliardo}. Let us fix $z \in \mathbb{R}^{m}$, $\varepsilon>0$ and let $v_\varepsilon \in \mathcal{A}_{\varepsilon,T,z}(Q_{R_{\varepsilon}};\mathbb{R}^{m})$ be such that 
    \begin{equation*}
        \mathcal{F}^{T}_{\varepsilon}(v_{\varepsilon},Q_{R_{\varepsilon}})=\varphi_{\varepsilon,T,R_{\varepsilon}}(z).
    \end{equation*}
    
    \noindent For every fixed $R>0$ and $\varepsilon$ small enough, by (\ref{26}) and by the fact that $\mathcal{F}^{T}_{\varepsilon}(v_{\varepsilon},\cdot)$ is an increasing set function, we have
    \begin{equation}
        \mathcal{F}^{T}_{\varepsilon}(v_{\varepsilon},Q_{R})\leq \mathcal{F}^{T}_{\varepsilon}(v_{\varepsilon},Q_{R_{\varepsilon}})\leq c_{2}|z|^{p}.
        \label{energia limitata}
    \end{equation}
    Now, we show that for $\varepsilon$ small enough
    \begin{equation}
    \label{ulteriore controllo}
        \mathcal{G}^{1,p}_{\varepsilon}(v_{\varepsilon},Q_{R})\leq CG^{1,p}_{\varepsilon}(v_{\varepsilon},Q_{R_{\varepsilon}})\le C |z|^{p}
    \end{equation}
    where $\mathcal{G}^{1,p}_{\varepsilon}$ is defined by \eqref{Gtondo} and the last inequality in \eqref{ulteriore controllo} is a consequence of Remark \ref{100}. The first inequality in  \eqref{ulteriore controllo} follows from the following argument. Let us define 
    $I_{\varepsilon}(Q_{R}):=\{\alpha \in \varepsilon\mathbb{Z}^{d}: (\alpha+[0,\varepsilon)^{d})\cap Q_{R}\neq \emptyset\}$ and $S^{k}_{\varepsilon,\alpha}$ as in \eqref{interfaccia}. It is clear that for $\varepsilon$ small enough
    \begin{equation*}
        Q_{R_{\varepsilon}}\supset\bigcup_{\alpha \in I_{\varepsilon}(Q_{R})}\alpha+[0,\varepsilon)^{d}\supset Q_{R},
    \end{equation*}
    and then, arguing as in the proof of Theorem \ref{Disuguaglianza di Gagliardo}, we get
    \begin{equation*}
    \begin{split}
    \mathcal{G}^{1,p}_{\varepsilon}(v_{\varepsilon},Q_{R})&=\int_{B_{1}}\int_{(Q_{R})_{\varepsilon}(\xi)}\left|\frac{v_{\varepsilon}(x+\varepsilon\xi)-v_{\varepsilon}(x)}{\varepsilon}\right|^{p}\,dx\,d\xi\\
    &\leq\int_{B_{1}}\int_{Q_{R}}\left|\frac{v_{\varepsilon}(x+\varepsilon\xi)-v_{\varepsilon}(x)}{\varepsilon}\right|^{p}\,dx\,d\xi\\
    &\leq\int_{B_{1}}\sum_{\alpha \in I_{\varepsilon}(Q_{R})}\int_{\alpha+[0,\varepsilon)^{d}}\left|\frac{v_{\varepsilon}(x+\varepsilon\xi)-v_{\varepsilon}(x)}{\varepsilon}\right|^{p}\,dx\,d\xi\\
    &\leq\int_{B_{1}}\sum_{k=1}^{d}\sum_{\alpha \in I_{\varepsilon}(Q_{R})}\int_{S^{k}_{\varepsilon,\alpha}}\left|\frac{v_{\varepsilon}(x+\varepsilon\xi)-v_{\varepsilon}(x)}{\varepsilon}\right|^{p}\,dx\,d\xi\\
    &\leq C\int_{B_{1}}|\xi|\, d\xi\sum_{k=1}^{d}\sum_{\alpha \in I_{\varepsilon}(Q_{R})}\varepsilon^{d}|D^{e_{k}}_{\varepsilon}v_{\varepsilon}(\alpha)|^{p}\\
    &\leq C G^{1,p}_{\varepsilon}(v_{\varepsilon},Q_{R_{\varepsilon}}).
    \end{split} 
    \end{equation*}
    Moreover, by Lemma \ref{53} and  Remark \ref{OSSPW}  applied with $E=Q_{1}$ and $\Omega=Q_{R}$, 
    we get
    \begin{equation}
        \label{equilimitatezza}
        \sup_{\varepsilon}\|v_{\varepsilon}\|_{L^{p}(Q_{R};\mathbb{R}^{m})}<+\infty.
    \end{equation}
    By \eqref{ulteriore controllo} and \eqref{equilimitatezza} we get 
    \begin{equation*}
        \sup_{\varepsilon}\{\|v_{\varepsilon}\|_{L^{p}(Q_{R};\mathbb{R}^{m})}+\mathcal{G}^{1,p}_{\varepsilon}(v_{\varepsilon},Q_{R})\}<+\infty.
    \end{equation*}
    and then, by Theorem \ref{compactness for concolution}, there exists a subsequence (not relabeled) $v_{\varepsilon}$ converging strongly in $L^{p}$ to $v \in W^{1,p}(Q_{R};\mathbb{R}^{m})$. By the arbitrariness of $R>0$ we infer that up to a subsequence $v_\varepsilon$ converges strongly in $L^{p}_{\loc}(\mathbb{R}^{d};\mathbb{R}^{m})$ to $v \in W^{1,p}_{\loc}(\mathbb{R}^{d};\mathbb{R}^{m})$  with $v \equiv 0$ in $Q_1$. As a direct consequence of Theorem \ref{Disuguaglianza di Gagliardo} we infer that $u:=v-z \in L^{p^{*}}(\mathbb{R}^{d};\mathbb{R}^{m})$. Indeed, let us define
    \begin{equation*}
            u_{\varepsilon}=
            \begin{cases}
                v_{\varepsilon}-z &\text{on}\ Q_{R_{\varepsilon}}\\
                0 &\text{on}\ \mathbb{R}^{d}\setminus Q_{R_{\varepsilon}}.
            \end{cases}
        \end{equation*}
        By (G0) and (\ref{energia limitata})  we get
        \begin{equation*}
            G^{1,p}_{\varepsilon}(u_{\varepsilon},\mathbb{R}^{d})=G^{1,p}_{\varepsilon}(u_{\varepsilon},Q_{R_{\varepsilon}})
                \leq C \mathcal{F}^{T}_{\varepsilon}(v_{\varepsilon},Q_{R_{\varepsilon}})\leq C|z|^{p}.
        \end{equation*}
    Then by \eqref{44}
\begin{equation*}
    \int_{\mathbb{R}^{d}}|u_{\varepsilon}|^{p^{*}}\leq CG^{1,p}_{\varepsilon}(u_{\varepsilon},\mathbb{R}^{d})\leq C|z|^{p},
\end{equation*}
and then by Fatou's lemma 
\begin{equation*}
    \int_{\mathbb{R}^{d}}|u|^{p^{*}}\leq \liminf_{\varepsilon \rightarrow 0}\int_{\mathbb{R}^{d}}|u_{\varepsilon}|^{p^{*}}<+\infty.
\end{equation*}

        \noindent By Theorem \ref{51} we deduce that 
        \begin{equation*}
            \liminf_{\varepsilon \rightarrow 0} \varphi_{\varepsilon,T,R_{\varepsilon}}(z)=\liminf_{\varepsilon \rightarrow 0} \mathcal{F}^{T}_{\varepsilon}(v_{\varepsilon},Q_{R_{\varepsilon}})
                \geq \liminf_{\varepsilon \rightarrow 0} \mathcal{F}^{T}_{\varepsilon}(v_{\varepsilon},Q_{R})\geq \int_{Q_{R}}f^{T}_{\hom}(\nabla u)\, dx.
        \end{equation*}
        Letting $R \rightarrow +\infty$ we obtain
        \begin{equation*}
            \liminf_{\varepsilon \rightarrow 0} \varphi_{\varepsilon,T,R_{\varepsilon}}(z)\geq \int_{\mathbb{R}^{d}}f^{T}_{\hom}(\nabla u)\, dx \geq \varphi^{T}(z).
        \end{equation*}
        We now prove 
        \begin{equation*}
            \limsup_{\varepsilon \rightarrow 0}\varphi_{\varepsilon,T,R_{\varepsilon}}(z)\leq \varphi^{T}(z).
        \end{equation*}
        A step of the proof of \cite[Proposition 6.7]{AGL}, shows that 
        \begin{equation*}
            \varphi^{T}(z)=\inf\bigg\{\int_{\mathbb{R}^{d}}f^{T}_{\hom}(\nabla u)\, dx: u\equiv -z\ \text{in}\ Q_{1}, u \in W^{1,p}(\mathbb{R}^{d};\mathbb{R}^{m})\ \text{compactly supported} \bigg\}.
        \end{equation*}
        Then, using a density argument, given $\eta >0$, there exists $u \in C^{\infty}_{c}(\mathbb{R}^{d};\mathbb{R}^{m})$ such that $u \equiv -z$ on $Q_{1}$ and 
        \begin{equation*}
            \int_{\mathbb{R}^{d}}f^{T}_{\hom}(\nabla u)\, dx \leq \varphi^{T}(z) + \eta.
        \end{equation*}
        Fix $R>0$ such that $\supp u \subset Q_{R}$. Then by \cite[Theorem 4.5]{alicic},  applied with $A=Q_{R}\setminus Q_{1}$, there exists a family of functions $u_{\varepsilon} \in L^{p}(Q_{R};\mathbb{R}^{m})$, with $u_{\varepsilon}\equiv -z$ on $Q_{1}$ and $u_{\varepsilon}\equiv 0$ on $\partial^{\varepsilon T}Q_{R}$ such that 
        \begin{equation*}
            \lim_{\varepsilon \rightarrow 0}\mathcal{F}^{T}_{\varepsilon}(u_{\varepsilon},\mathbb{R}^{d})=\int_{\mathbb{R}^{d}}f^{T}_{\hom}(\nabla u)\, dx.
        \end{equation*}
        Hence, 
        \begin{equation*}
            \limsup_{\varepsilon \rightarrow 0}\varphi_{\varepsilon,T,R_{\varepsilon}}(z)= \limsup_{\varepsilon \rightarrow 0}\mathcal{F}^{T}_{\varepsilon}(v_{\varepsilon},Q_{R_{\varepsilon}})
            \leq \int_{\mathbb{R}^{d}}f^{T}_{\hom}(\nabla u)\, dx
            \leq  \varphi^{T}(z)+\eta.
        \end{equation*}
        By the arbitrariness of $\eta >0$, we get
        \begin{equation*}
            \limsup_{\varepsilon \rightarrow 0}\varphi_{\varepsilon,T,R_{\varepsilon}}(z)\leq \varphi^{T}(z),
        \end{equation*}
        and then the thesis.
    \end{proof}
\label{DC}
\end{proposition}

\begin{proposition}
    Let $\varepsilon_{j}\rightarrow 0$ and $R_{j}\rightarrow +\infty$ as $j \rightarrow +\infty$ and let $\{u_{j}\}_{j}\subset \mathcal{A}_{\varepsilon_{j}}(\Omega;\mathbb{R}^{m})$ be a bounded sequence in $L^{\infty}(\Omega;\mathbb{R}^{m})$ such that $\sup_{j}\mathcal{F}^{T}_{\varepsilon
    _{j}}(u_{j})<+\infty$. Assume that $u_{j}\rightarrow u$ in $L^{p}(\Omega;\mathbb{R}^{m})$ for some $u \in W^{1,p}(\Omega;\mathbb{R}^{m})$. Let $\{\rho_{j}\}_{j}$ be a sequence of the form $\rho_{j}=\beta\delta_{j}$, with $\beta<\frac{1}{2}$ and for all $i \in Z_{j}(\Omega)$, where $Z_{j}(\Omega)$ is defined in \eqref{insieme di indici}, let $u^{i}_{j}$ be defined as in Lemma \ref{JL}. Set $Q^{i}_{j}:=i\delta_{j}+\bigg(\bigg[-\frac{\delta_{j}}{2},\frac{\delta_{j}}{2}\bigg)^{d}\cap \varepsilon_{j}\mathbb{Z}^{d}\bigg)$ and let $\psi^{T}_{j}\in \mathcal{A}_{\varepsilon_{j}}(\Omega;\mathbb{R}^{m})$ be defined for $\alpha \in \Omega_{j}$ as 
    \begin{equation*}
        \psi^{T}_{j}(\alpha)=\sum_{i\in Z_{j}(\Omega)}\varphi_{\varepsilon_{j},T,R_{j}}(u^{i}_{j})\chi_{Q^{i}_{j}}(\alpha).
    \end{equation*}
     Then,
    \begin{equation*}
        \lim_{j\rightarrow +\infty} \sum_{\alpha \in \Omega_{j}}\varepsilon_{j}^{d}\psi^{T}_{j}(\alpha)=\int_{\Omega}\varphi^{T}(u)\, dx.  
     \end{equation*}

    \begin{proof}
        First of all we observe that
        \begin{equation*}
            \sum_{\alpha \in \Omega_{j}}\varepsilon_{j}^{d}\psi^{T}_{j}(\alpha)=\sum_{i \in Z_{j}(\Omega)}\delta^{d}_{j}\varphi_{\varepsilon_{j},T,R_{j}}(u^{i}_{j}).
        \end{equation*}
        Indeed,
        \begin{equation*}
        \begin{split}
             \sum_{\alpha \in \Omega_{j}}\varepsilon_{j}^{d}\psi^{T}_{j}(\alpha)&=\sum_{\alpha \in \Omega_{j}}\sum_{i\in Z_{j}(\Omega)}\varepsilon_{j}^{d}\varphi_{\varepsilon_{j},T,R_{j}}(u^{i}_{j})\chi_{Q^{i}_{j}}(\alpha)\\
             &=\sum_{i \in Z_{j}(\Omega)}\varepsilon_{j}^{d}\bigg(\frac{\delta_{j}}{\varepsilon_{j}}\bigg)^{d}\varphi_{\varepsilon_{j},T,R_{j}}(u^{i}_{j})\\
             &=\sum_{i \in Z_{j}(\Omega)}\delta_{j}^{d}\varphi_{\varepsilon_{j},T,R_{j}}(u^{i}_{j}).
             \end{split}
        \end{equation*}
            Now, we get the following estimate
            \begin{equation*}
    \begin{split}
        \bigg| \sum_{i \in Z_{j}(\Omega)} &\delta_{j}^{d} \varphi_{\varepsilon_{j},T,R_{j}}(u^{i}_{j}) 
        - \int_{\Omega} \varphi^{T}(u(x))\, dx \bigg| \\
        &=\bigg| \sum_{i \in Z_{j}(\Omega)} \varphi_{\varepsilon_{j},T,R_{j}}(u^{i}_{j}) \int_{Q^{i}_{j}}\, dx 
        - \sum_{i \in Z_{j}(\Omega)} \int_{Q^{i}_{j}} \varphi^{T}(u(x))\, dx -\int_{\Omega \setminus \cup_{i \in Z_{j}(\Omega)}Q^{i}_{j}}\varphi^{T}(u(x))\, dx\bigg| \\
        &\leq \sum_{i \in Z_{j}(\Omega)} \int_{Q^{i}_{j}} \left| \varphi_{\varepsilon_{j},T,R_{j}}(u^{i}_{j}) 
        - \varphi^{T}(u(x)) \right|\, dx + \int_{\Omega \setminus \cup_{i \in Z_{j}(\Omega)}Q^{i}_{j}}|\varphi^{T}(u(x))|\, dx\\
        &\leq \sum_{i \in Z_{j}(\Omega)} \sum_{\alpha \in Q^{i}_{j}} \int_{\alpha + [0,\varepsilon_j)^d} 
        \left| \varphi_{\varepsilon_{j},T,R_{j}}(u^{i}_{j}) - \varphi_{\varepsilon_{j},T,R_{j}}(u_{j}(x)) \right|\, dx \\
        &\ + \sum_{i \in Z_{j}(\Omega)} \sum_{\alpha \in Q^{i}_{j}} \int_{\alpha + [0,\varepsilon_j)^d} 
        \left| \varphi_{\varepsilon_{j},T,R_{j}}(u_{j}(x)) - \varphi^{T}(u(x)) \right|\, dx+\int_{\Omega \setminus \cup_{i \in Z_{j}(\Omega)}Q^{i}_{j}}|\varphi^{T}(u(x))|\, dx \\
        &\leq \sum_{i \in Z_{j}(\Omega)} \sum_{\alpha \in Q^{i}_{j}} \varepsilon_{j}^d \left| \varphi_{\varepsilon_{j},T,R_{j}}(u^{i}_{j}) 
        - \varphi_{\varepsilon_{j},T,R_{j}}(u_{j}(\alpha)) \right|  \\
        & + \sum_{i \in Z_{j}(\Omega)} \sum_{\alpha \in Q^{i}_{j}} \varepsilon_{j}^d \left| \varphi_{\varepsilon_{j},T,R_{j}}(u_{j}(\alpha)) 
        - \varphi^{T}(u_{j}(\alpha)) \right|\\
        &+  \sum_{i \in Z_{j}(\Omega)} \sum_{\alpha \in Q^{i}_{j}} \int_{\alpha + [0,\varepsilon_j)^d} \left|\varphi^{T}(u_{j}(x))-\varphi^{T}(u(x))\right|\, dx + \int_{\Omega \setminus \cup_{i \in Z_{j}(\Omega)}Q^{i}_{j}}|\varphi^{T}(u(x))|\, dx \\
        &:= S_{j}^{1} + S_{j}^{2}+S_{j}^{3}+S_{j}^{4}.
    \end{split}
\end{equation*}

            \noindent By hypothesis  $u_{j} \in L^{\infty}(\Omega;\mathbb{R}^{m})$, then by Proposition \ref{DC} we deduce that $S^{2}_{j}\rightarrow 0$. Since $\left|\Omega \setminus \cup_{i \in Z_{j}(\Omega)}Q^{i}_{j}\right|\rightarrow 0$ we also infer that $S_{j}^{4}\rightarrow 0$. Moreover, by the fact that $\varphi^{T}$ is locally Lipschitz and $u_{j}\rightarrow u$ in $L^{p}$ we get that $S_{j}^{3}\rightarrow 0$. Finally, by Proposition \ref{200}, we may estimate $S^{1}_{j}$ as follows:
            \begin{equation*}
                S^{1}_{j}\leq C \sum_{i \in Z_{j}(\Omega)}\sum_{\alpha \in Q^{i}_{j}}\varepsilon_{j}^{d}|u^{i}_{j}-u_{j}(\alpha)|.
            \end{equation*}
            Using a discrete version of Hölder inequality, we get
            \begin{equation*}
                \begin{split}
                    \sum_{i \in Z_{j}(\Omega)}\sum_{\alpha \in Q^{i}_{j}}\varepsilon_{j}^{d}|u^{i}_{j}-u_{j}(\alpha)|&\leq \varepsilon_{j}^{d}\sum_{i \in Z_{j}(\Omega)} (\sharp Q^{i}_{j})^{1-\frac{1}{p}}\bigg(\sum_{\alpha \in Q^{i}_{j}}|u^{i}_{j}-u_{j}(\alpha)|^{p}\bigg)^{\frac{1}{p}}\\
                    &=\delta_{j}^{d}\delta_{j}^{-\frac{d}{p}}\sum_{i \in Z_{j}(\Omega)}\bigg(\sum_{\alpha \in Q^{i}_{j}}\varepsilon_{j}^{d}|u^{i}_{j}-u_{j}(\alpha)|^{p}\bigg)^{\frac{1}{p}}.
                \end{split}
            \end{equation*}
            By Lemma \ref{53bis} and Remark \ref{OSSPW} we deduce that 
            \begin{equation*}
                \sum_{\alpha \in Q^{i}_{j}}\varepsilon_{j}^{d}|u^{i}_{j}-u_{j}(\alpha)|^{p}\leq C\delta_{j}^{p}\sum_{k=1}^{d}\sum_{\alpha \in (Q^{i}_{j})_{\varepsilon_{j}}(e_{k})}\varepsilon_{j}^{d}\bigg|\frac{u_{j}(\alpha+e_{k}\varepsilon_{j})-u_{j}(\alpha)}{\varepsilon_{j}}\bigg|^{p}.
            \end{equation*}
            Then using the fact that the function $x \mapsto x^{\frac{1}{p}}$ is concave, we get
            \begin{equation*}
            \begin{split}
                \sum_{i \in Z_{j}(\Omega)}\sum_{\alpha \in Q^{i}_{j}}&\varepsilon_{j}^{d}|u^{i}_{j}-u_{j}(\alpha)|\leq C\delta_{j}^{d}\delta_{j}^{-\frac{d}{p}}\sum_{i \in Z_{j}(\Omega)}\delta_{j}\bigg(\sum_{k=1}^{d}\sum_{\alpha\in (Q^{i}_{j})_{\varepsilon_{j}}(e_{k})}\varepsilon_{j}^{d}\bigg|\frac{u_{j}(\alpha+e_{k}\varepsilon_{j})-u_{j}(\alpha)}{\varepsilon_{j}}\bigg|^{p}\bigg)^{\frac{1}{p}}\\
                &\leq C\delta_{j}^{d}\delta_{j}^{-\frac{d}{p}}\delta_{j}(\sharp Z_{j}(\Omega))\bigg(\frac{1}{(\sharp Z_{j}(\Omega))}\sum_{i \in Z_{j}(\Omega)}\sum_{k=1}^{d}\sum_{\alpha \in (Q^{i}_{j})_{\varepsilon_{j}}(e_{k})}\varepsilon_{j}^{d}\bigg|\frac{u_{j}(\alpha+e_{k}\varepsilon_{j})-u_{j}(\alpha)}{\varepsilon_{j}}\bigg|^{p}\bigg)^{\frac{1}{p}}\\
                &= C\delta_{j}^{d}\delta_{j}^{-\frac{d}{p}}\delta_{j}(\sharp Z_{j}(\Omega))^{1-\frac{1}{p}}\bigg(\sum_{i \in Z_{j}(\Omega)}\sum_{k=1}^{d}\sum_{\alpha \in (Q^{i}_{j})_{\varepsilon_{j}}(e_{k})}\varepsilon_{j}^{d}\bigg|\frac{u_{j}(\alpha+e_{k}\varepsilon_{j})-u_{j}(\alpha)}{\varepsilon_{j}}\bigg|^{p}\bigg)^{\frac{1}{p}}\\
                &\leq C\delta_{j}(G^{1,p}_{\varepsilon_{j}}(u_{j},\Omega))^{\frac{1}{p}}\leq C \delta_{j}(\mathcal{F}^{T}_{\varepsilon_{j}}(u_{j}))^{\frac{1}{p}}\leq C \delta_{j}(\sup_{j}\mathcal{F}^{T}_{\varepsilon_{j}}(u_{j}))^{\frac{1}{p}} \leq C\delta_{j}.
            \end{split}
            \end{equation*}
            Hence,
            \begin{equation*}
                S^{1}_{j}\leq C\delta_{j}\rightarrow 0.
            \end{equation*}
    \end{proof}
    \label{lastool}
\end{proposition}
\section{Proof of Theorem \ref{MAIN}}

 We will divide the proof of Theorem \ref{MAIN} into three steps. We first show that it suffices to prove the theorem for the truncated functionals $F^{T}_{\varepsilon,\delta_{\varepsilon}}$, and then we will deal separately with the $\Gamma$-$\liminf$ and the $\Gamma$-$\limsup$ inequalities.\\

\noindent \textbf{Step 1}\\
It is not restrictive to prove the theorem under the additional assumption that there exists $T>2$ such that $f(\xi,z)=0$ if $|\xi|>T$. Indeed, under the hypothesis of Theorem \ref{MAIN}, assume that for every $T>2$
\begin{equation*}
    \Gamma(L^{p})-\lim_{\varepsilon}F^{T}_{\varepsilon,\delta_{\varepsilon}}(u)=
    \begin{cases}
    \displaystyle
        \int_{\Omega}f^{T}_{\hom}(\nabla u)\, dx+\gamma^{d-p} \int_{\Omega}\varphi^{T}(u)\, dx &\text{if}\ u \in W^{1,p}(\Omega;\mathbb{R}^{m})\\
        +\infty &\text{otherwise},
    \end{cases}
\end{equation*}
where $F^{T}_{\varepsilon,\delta_{\varepsilon}}$ is defined by \eqref{funzionale troncato} and $f^{T}_{\hom}$ and $\varphi^{T}$ are defined by \eqref{f omogenizzata} and \eqref{denista troncata}, respectively. If we prove that
\begin{equation}
    \Gamma(L^{p})-\lim_{\varepsilon}F_{\varepsilon,\delta_{\varepsilon}}(u)=\lim_{T\rightarrow +\infty} \Gamma(L^{p})-\lim_{\varepsilon}F^{T}_{\varepsilon,\delta_{\varepsilon}}(u).
    \label{T infinito}
\end{equation}
then, by Monotone Convergence Theorem, the statement follows once we prove that for every $M \in \mathbb{R}^{d \times m}$ and $z \in \mathbb{R}^{m}$
\begin{equation*}
    \lim_{T \rightarrow +\infty} f^{T}_{\hom}(M)=f_{\hom}(M),\ \quad \lim_{T \rightarrow +\infty}\varphi^{T}(z)=\varphi(z).
\end{equation*}
The inequalities above are again a straightforward consequence of Monotone Convergence Theorem. In order to prove (\ref{T infinito}) it suffices to show that
\begin{equation*}
        \Gamma(L^{p})-\limsup_{\varepsilon\rightarrow 0}F_{\varepsilon,\delta_{\varepsilon}}(u)\leq \lim_{T \rightarrow +\infty}\Gamma(L^{p})-\lim_{\varepsilon \rightarrow 0}F^{T}_{\varepsilon,\delta_{\varepsilon}}(u).
\end{equation*}

\noindent For reader's convenience we will specify the dependence of our energies on the open sets $A \subset \mathbb{R}^{d}$. Setting $F^{''}(u,\Omega):=\Gamma(L^{p})-\limsup_{\substack{\varepsilon \rightarrow 0}}F_{\varepsilon,\delta_{\varepsilon}}(u,\Omega)$ and $F^{T}(u,\Omega):=\Gamma(L^{p})-\lim_{\substack{\varepsilon \rightarrow 0}}F^{T}_{\varepsilon,\delta_{\varepsilon}}(u,\Omega)$, let us prove that for every $\eta>0$ there exists $T_{\eta}>0$ such that for every $u \in W^{1,p}(\Omega,\mathbb{R}^{m})$
\begin{equation*}
    F^{''}(u,\Omega)\leq F^{T_{\eta}}(u,\Omega)+C\eta.
\end{equation*}
By a density argument it is sufficient to prove the inequality for $u \in C^{\infty}_{c}(\mathbb{R}^{d},\mathbb{R}^{m})$. Let us fix $\eta>0$, then there exists $T_{\eta}>0$ such that $\sum_{|\xi|>T_{\eta}}M(\xi)<\eta$. By the external regularity of $F^{T}(u,\Omega)$ there exists $\Omega'\supset \supset \Omega$ such that
\begin{equation}
    \label{regolarità esterna}
    F^{T_{\eta}}(u,\Omega')\leq F^{T_{\eta}}(u,\Omega)+\eta.
\end{equation}
Consider $u_{\varepsilon}\rightarrow u$ in $L^{p}(\Omega';\mathbb{R}^{m})$ such that $F_{\varepsilon}^{T_{\eta}}(u_{\varepsilon},\Omega')$ converges to $F^{T_{\eta}}(u,\Omega')$. Observe that, from \eqref{crescita dell' omogenizzata},  $f^T_{\hom}(M)\le f_{\hom}(M)\le C|M|^p$ for every $M\in\mathbb{R}^{m\times d}$. Moreover by Proposition \ref{non degeneracy of capacitary density} and Proposition \ref{DC}, $\varphi^T(z)\le c_2|z|^p$, for every $z\in\R^m$. Hence   
\begin{equation}
    F_{\varepsilon}^{T_{\eta}}(u_{\varepsilon},\Omega')\leq C,
    \label{equilimitatezza delle troncate}
\end{equation}
with $C$ independent of $\eta$. We split our energies as follows

   \begin{equation*}
       F_{\varepsilon}(u_{\varepsilon},\Omega)=F_{\varepsilon}^{T_{\eta}}(u_{\varepsilon},\Omega)+F_{\varepsilon}^{>T_{\eta}}(u_{\varepsilon},\Omega)
   \end{equation*}

   \noindent where 
   \begin{equation*}
   F_{\varepsilon}^{>T_{\eta}}(u_{\varepsilon},\Omega):=\sum_{|\xi|>T_{\eta}}\sum_{\alpha \in \Omega_{\varepsilon}(\xi) }\varepsilon^{d}f(\xi,D^{\xi}_{\varepsilon}u(\alpha)).
   \end{equation*}
   Let us define $\tilde{\Omega}'_{\varepsilon}:=\{x \in \Omega': \dist(x,\partial \Omega')>2 \sqrt{d}\varepsilon\}$. If $\varepsilon$ is sufficiently small we have $\Omega \subset \subset \tilde{\Omega}'_{\varepsilon}$. By Lemma \ref{AC} applied with $\tilde{\Omega}'_{\varepsilon}$, assumption (G) and \eqref{equilimitatezza delle troncate} we get
   \begin{equation}
       \label{stima fuori dalle troncate}
       F_{\varepsilon}^{>T_{\eta}}(u_{\varepsilon},\Omega)\leq \sum_{|\xi|>T}M(\xi)G^{1,p}_{\varepsilon}(u_{\varepsilon},\Omega')\leq C \eta.
       \end{equation}
      Hence, by \eqref{regolarità esterna} and \eqref{stima fuori dalle troncate}

      \begin{equation*}
       \begin{split}
          F^{''}(u,\Omega)&\leq\limsup_{\varepsilon \rightarrow 0} F_{\varepsilon}(u_{\varepsilon},\Omega)\leq \limsup_{\varepsilon \rightarrow 0} F^{T_{\eta}}_{\varepsilon}(u_{\varepsilon},\Omega')+C\eta \\
          &= F^{T_{\eta}}(u,\Omega')+C\eta \leq F^{T_{\eta}}(u,\Omega)+C\eta.
        \end{split}
      \end{equation*}

 \noindent \textbf{Step 2}. Now we prove the validity of the $\Gamma$-$\liminf$ inequality for $F^{T}_{\varepsilon,\delta_{\varepsilon}}$.\\
\noindent Given $\varepsilon_{j}\rightarrow 0$ as $j \rightarrow +\infty$, let $u \in W^{1,p}(\Omega;\mathbb{R}^{m})$ and $u_j \rightarrow u$ in $L^{p}(\Omega;\mathbb{R}^{m})$ be such that $\sup_{j}F^{T}_{\varepsilon_{j},\delta_{j}}(u_j)<+\infty$. We first assume that $\{u_{j}\}_j$ is bounded in $L^{\infty}$ and then we will remove this assumption using Lemma \ref{TL}. \\
We fix $k \in \mathbb{N}$ and we consider a sequence $\{\rho_{j}\}_j$ of the form $\rho_{j}=\beta\delta_{j}$ with $\beta<\frac{1}{2}$. Let $Z_{j}(\Omega)$ and $\rho^{i}_{j}$ be defined as in Lemma \ref{JL}. We apply Lemma \ref{JL} to $\{u_{j}\}_j$ in order to get a new sequence $\{w_j\}_j$ such that $w_{j} \rightarrow u$ in $L^{p}(\Omega;\mathbb{R}^{m})$ and it satisfies (\ref{outcornici})-(\ref{diffenergia}). In particular the following inequality holds
\begin{equation}
    \liminf_{j\rightarrow +\infty}F^{T}_{\varepsilon_{j},\delta_j}(u_{j})\geq  \liminf_{j\rightarrow +\infty}F^{T}_{\varepsilon_{j},\delta_j}(w_{j})-\frac{C}{k}.
    \label{disuguaglianza JL}
\end{equation}
Set
        \begin{equation*}
            E_{j}=\bigcup_{i \in Z_{j}(\Omega)}Q(i\delta_{j},\rho^{i}_{j}).
        \end{equation*}
        Note that we have 
        \begin{equation*}
            \label{splitted energy}\liminf_{j\rightarrow +\infty}F^{T}_{\varepsilon_{j},\delta_j}(w_{j})\geq \liminf_{j\rightarrow +\infty}F^{T}_{\varepsilon_j,\delta_j}(w_{j},E_{j})+\liminf_{j \rightarrow +\infty} F^{T}_{\varepsilon_{j},\delta_j}(w_{j},\Omega\setminus E_{j}).
        \end{equation*}
        We first estimate the energy contribution of $\{u_{j}\}_{j}$ on $\Omega\setminus E_{j}$, i.e. far from the perforations. To this end we use the argument exploited in \cite[Section 7, Subsection A.1]{Siga}. Let us define $v_{j}\in {A}_{\varepsilon_{j}}(\Omega;\mathbb{R}^{m})$ by modifying $w_{j}$ as follows
        \begin{equation*}
            v_{j}(\alpha)=
            \begin{cases}
                u^{i}_{j} &\alpha \in Q(i\delta_{j},\rho^{i}_{j}),\ i \in Z_{j}(\Omega)\\
                w_{j}(\alpha) &\text{otherwise}.
            \end{cases}
        \end{equation*}
        Note that since $u_j \in L^{\infty}(\Omega;\mathbb{R}^{m})$, by construction also $v_{j}, w_{j} \in L^{\infty}(\Omega;\mathbb{R}^{m})$. Now we prove that $v_{j}\rightarrow u$ in $L^{p}(\Omega;\mathbb{R}^{m})$. By Lemma \ref{53bis} and using the fact that $u_{j}\rightarrow u$ in $L^{p}(\Omega;\mathbb{R}^{m})$, $w_{j}\rightarrow u$ in $L^{p}(\Omega;\mathbb{R}^{m})$, we have
        \begin{equation*}
            \begin{split}
                \int_{\Omega}|v_{j}-u_{j}|^{p}\, dx&=\sum_{\alpha \in \Omega_{j}}\varepsilon_{j}^{d}|u_{j}(\alpha)-v_{j}(\alpha)|^{p}+o(1)\\
                &\leq \bigg[\sum_{\alpha \in \Omega_{j}\setminus E_{j}}\varepsilon_{j}^{d}|u_{j}(\alpha)-v_{j}(\alpha)|^{p}+\sum_{\alpha\in E_{j}}\varepsilon_{j}^{d}|u_{j}(\alpha)-v_{j}(\alpha)|^{p}\bigg]+o(1)\\
                &\leq \bigg[\sum_{\alpha \in \Omega_{j}\setminus E_{j}}\varepsilon_{j}^{d}|u_{j}(\alpha)-w_{j}(\alpha)|^{p}+\sum_{\alpha\in E_{j}}\varepsilon_{j}^{d}|u_{j}(\alpha)-v_{j}(\alpha)|^{p}\bigg]+o(1)\\
                &\leq \bigg[\int_{\Omega}|u_{j}-w_j|^{p}\,dx+\sum_{i \in Z_{j}(\Omega)}\sum_{\alpha \in Q(i\delta_{j},\rho^{i}_{j})}\varepsilon_{j}^{d}|u_{j}(\alpha)-u^{i}_{j}|^{p}\bigg]+o(1)\\
                &\leq C\delta_{j}^{p}\sum_{i\in Z_{j}(\Omega)}\sum_{k=1}^{d}\sum_{\alpha \in (Q(i\delta_{j},\rho^{i}_{j}))_{\varepsilon_{j}}(e_{k})}\varepsilon_{j}^{d}\bigg|\frac{u_{j}(\alpha +e_{k}\varepsilon_{j})-u_{j}(\alpha)}{\varepsilon_{j}}\bigg|^{p}+o(1)\\
&=C\delta_{j}^{p}G^{1,p}_{\varepsilon_{j}}(u_{j},\Omega)+o(1)
                \leq C\delta_{j}^{p}F^{T}_{\varepsilon_{j},\delta_j}(u_{j})+o(1)
                \leq C\delta_{j}^{p}+o(1)
            \end{split}
        \end{equation*}
        By construction
        \begin{equation*}
        F^{T}_{\varepsilon_j,\delta_j}(w_j,\Omega_j\setminus E_j)\geq F^{T}_{\varepsilon_j,\delta_j}(v_j).
        \end{equation*}
        Then by Theorem \ref{51}, we conclude that
        \begin{equation}
    \liminf_{j}F^{T}_{\varepsilon_{j},\delta_j}(w_{j};\Omega_{j}\setminus E_{j})\geq \liminf_{j}F^{T}_{\varepsilon_{j},\delta_j}(v_{j})
            \geq \int_{\Omega}f^T_{\hom}(\nabla u)\, dx. 
            \label{estimate far from the ps}
        \end{equation}
         
        \noindent Now we focus our attention on the contribution of $\{u_{j}\}_j$ on $E_{j}$, i.e. close to the perforations. Set $R^{i}_{j}:=\frac{\rho^{i}_{j}}{r_{\delta_{j}}}$ and $s_{j}:=\frac{\varepsilon_{j}}{r_{\delta_{j}}}$. Since $\varepsilon_j=o(r_{\delta_{j}})$, $s_{j}\rightarrow 0$ when $j \rightarrow +\infty$. For  $i \in Z_{j}(\Omega)$ we define the functions $v_{i,j}: s_{j}\mathbb{Z}^{d}\rightarrow \mathbb{R}^{m}$ as follows
        
        \begin{equation*}
    v_{i,j}(\alpha) =
    \begin{cases}
        w_{j}(i\delta_{j} + r_{\delta_{j}}\alpha) & \text{if } \alpha \in Q_{R^{i}_{j}} \\
        u^{i}_{j} & \text{if } \alpha \in s_{j}\mathbb{Z}^{d} \setminus Q_{R^{i}_{j}}.
    \end{cases}
\end{equation*}
Using hypothesis (H) we get
\begin{equation*}
    \begin{split}
        F^{T}_{\varepsilon_{j},\delta_j}(w_{j},Q(i\delta_{j},\rho^{i}_{j}))&=\sum_{|\xi|\leq T}\sum_{\alpha \in (Q(i\delta_{j},\rho^{i}_{j}))_{\varepsilon_{j}}(\xi)}\varepsilon_{j}^{d}f\bigg(\xi,\frac{w_{j}(\alpha + \varepsilon_{j}\xi)-w_{j}(\alpha)}{|\xi|\varepsilon_j}\bigg)\\
        &=\sum_{|\xi|\leq T}\sum_{\alpha \in \big(Q_{R^{i}_{j}}\big)_{s_{j}}(\xi)}\varepsilon_{j}^{d}f\bigg(\xi,\frac{v_{i,j}(\alpha+s_j\xi)-v_{i,j}(\alpha)}{r_{\delta_j}s_j|\xi|}\bigg)\\
        &=\sum_{|\xi|\leq T}\sum_{\alpha \in \big(Q_{R^{i}_{j}}\big)_{s_{j}}(\xi)}\frac{\varepsilon_{j}^{d}}{r^{p}_{\delta_j}}f\bigg(\xi,\frac{v_{i,j}(\alpha+s_j\xi)-v_{i,j}(\alpha)}{s_j|\xi|}\bigg)\\
        &=\sum_{|\xi|\leq T}\sum_{\alpha \in \big(Q_{R^{i}_{j}}\big)_{s_{j}}(\xi)}r^{d-p}_{\delta_j}s^{d}_{j}f\bigg(\xi,\frac{v_{i,j}(\alpha+s_j\xi)-v_{i,j}(\alpha)}{s_j|\xi|}\bigg)\\
        &=r_{\delta_j}^{d-p}F^{T}_{s_j,\delta_j}(v_{i,j},Q_{R^{i}_{j}})\\
        &\geq r_{\delta_j}^{d-p}\varphi_{s_j,T,R^{i}_{j}}(u^{i}_{j}).
    \end{split}
\end{equation*}
    Summing over $i \in Z_j(\Omega)$,
    \begin{equation*}
        \sum_{i \in Z_j(\Omega)}F^{T}_{s_j,\delta_j}(v_{i,j},Q_{R^{i}_{j}})\geq\sum_{i \in Z_j(\Omega)}r_{\delta_j}^{d-p}\varphi_{s_j,T,R^{i}_{j}}(u^{i}_{j}).
    \end{equation*}
    Note that the right hand side above, can be written as
    \begin{equation*}
        r_{\delta_j}^{d-p}\sum_{i \in Z_j(\Omega)}\varphi_{s_j,T,R^{i}_{j}}(u^{i}_{j})=\frac{r_{\delta_j}^{d-p}}{\bigg(\frac{\delta_j}{\varepsilon_j}\bigg)^{d}}\sum_{\alpha \in \Omega_j}\sum_{i \in Z_j(\Omega)}\varphi_{s_j,T,R^{i}_{j}}(u^{i}_{j})\chi_{Q^{i}_{j}}(\alpha)
        =\Bigg(\frac{r_{\delta_j}}{\delta_j^{\frac{d}{d-p}}}\Bigg)^{d-p}\sum_{\alpha \in \Omega_j}\varepsilon_{j}^{d}\psi^{T}_{j}(\alpha).
    \end{equation*}
    By Proposition \ref{lastool} we deduce that 
    \begin{equation}
        \liminf_{j}F^{T}_{\varepsilon_j,\delta_j}(w_j,E_j)\geq \liminf_{j}\sum_{i \in Z_j(\Omega)}r_{\delta_j}^{d-p}\varphi_{s_j,T,R^{i}_{j}}(u^{i}_{j})\geq \gamma^{d-p}\int_{\Omega}\varphi^{T}(u(x))\, dx.
        \label{estimate on the ps}
    \end{equation}
    
    \noindent By (\ref{disuguaglianza JL}), (\ref{estimate far from the ps}), (\ref{estimate on the ps}) and by the arbitrariness of $k$, we finally get
    \begin{equation*}
        \liminf_{j}F^{T}_{\varepsilon_j,\delta_j}(u_j)\geq \int_{\Omega}f^T_{\hom}(\nabla u)+\gamma^{d-p}\int_{\Omega}\varphi^T(u)\, dx.
    \end{equation*}
    It remains to show that the $\Gamma$-liminf inequality holds even if we remove the boundedness assumption on the sequence $\{u_{j}\}$. For all $L \in \mathbb{N}$ and $\eta >0$ we apply the previous arguments to the truncated sequence $\{t_{L}(u_{j})\}_{j}$, where $t_{L}$ is as in the statement of Lemma \ref{TL} and in particular it is such that
    \begin{equation}
        \label{inequality TL}
        \liminf_{j}F^T_{\varepsilon_j,\delta_j}(t_L(u_j))\leq \liminf_{j}F^T_{\varepsilon_j,\delta_j}(u_j)+\eta.
    \end{equation}
    
    \noindent Then, we get
    \begin{equation*}
        \liminf_{j}F^{T}_{\varepsilon_j,\delta_{\varepsilon_j}}(t_L(u_j))\geq \int_{\Omega}f^T_{\hom}(\nabla t_L(u))\, dx+\gamma^{d-p}\int_{\Omega}\varphi^T(t_{L}(u))\, dx.
    \end{equation*}
    Note that $t_L(u) \rightarrow u$ as $L \rightarrow +\infty$, with respect to the strong convergence of $W^{1,p}(\Omega;\mathbb{R}^{m})$. By (\ref{inequality TL}) and the arbitrariness of $\eta$, we can pass to the limit as $L \rightarrow +\infty$ and we get 
    \begin{equation*}
        \liminf_{j}F^{T}_{\varepsilon_j,\delta_j}(u_{j})\geq F^{T}(u).
    \end{equation*}

\noindent \textbf{Step 3}. Finally we prove the validity of the $\Gamma$-limsup inequality.\\
\noindent By a density argument it is sufficient to prove this inequality for $u \in C^{\infty}_{c}(\mathbb{R}^{d};\mathbb{R}^{m})$. For such $u$ let us consider $\Omega'\in \mathcal{A}^{\reg}(\mathbb{R}^{d})$ such that $\Omega'\supset \supset \Omega$. By Theorem \ref{51} there exists $\tilde{u}_j \in \mathcal{A}_{\varepsilon_j}(\Omega',\mathbb{R}^m)$, $\tilde{u}_j \rightarrow u$ in $L^{p}(\Omega',\mathbb{R}^{m})$ such that
         \begin{equation}
             \lim_{j\rightarrow +\infty}\mathcal{F}^{T}_{\varepsilon_j}(\tilde{u}_j,\Omega')=\int_{\Omega'}f^{T}_{\hom}(\nabla u).
             \label{recovery sequence}
         \end{equation}
         By Lemma \ref{TL} we may suppose also that $\|\tilde{u}_j\|_{L^{\infty}(\Omega',\mathbb{R}^{m})}<+\infty$. Given $k \in \mathbb{N}$, by Lemma \ref{JL} there exists $\{w_{j}\}_j$ satisfying \eqref{outcornici}-\eqref{diffenergia} such that 
         \begin{equation}
             \mathcal{F}^{T}_{\varepsilon_{j}}(w_j,\Omega')\leq \mathcal{F}^{T}_{\varepsilon_{j}}(\tilde{u}_j,\Omega')+\frac{C}{k}.\label{inequality JL}
         \end{equation}
         Now we consider $E_{j}=\bigcup_{i \in Z_j(\Omega')}Q(i\delta_j,\rho^{i}_{j})$. As in Step 2, let us set $R^{i}_{j}:=\frac{\rho^{i}_{j}}{r_{\delta_j}}$ and $s_j:=\frac{\varepsilon_j}{r_{\delta_j}}$. For $i \in Z_j(\Omega')$, let $\tilde{v}_{j,i} \in \mathcal{A}_{s_j,T,\tilde{u}^{i}_{j}}(Q_{R^{i}_{j}};\mathbb{R}^{m})$ such that 
         \begin{equation}
             F^{T}_{s_j,\delta_j}(\tilde{v}_{j,i},Q_{R^{i}_{j}})=\varphi_{s_{j},T,R^{i}_{j}}(\tilde{u}^{i}_{j}).
             \label{minimizing sequence}
         \end{equation}
        
        \noindent Set
        \begin{equation*}
            v_{j,i}(\alpha)=\tilde{v}_{j,i}\bigg(\frac{\alpha-i\delta_j}{r_{\delta_{j}}}\bigg), \quad \alpha \in Q(i\delta_j,\rho^{i}_{j}),
        \end{equation*}
        and then
        \begin{equation*}
            u_{j}(\alpha)=
            \begin{cases}
                w_{j}(\alpha)&\text{if}\ \alpha \in \Omega'\setminus E_{j}\\
                v_{j,i}(\alpha) &\text{if}\ \alpha \in Q(i\delta_j,\rho^{i}_{j}).
            \end{cases}
        \end{equation*}
        By a change of variable we get
        \begin{equation}
            \label{cambio di variabili sulle perforazioni}F^{T}_{\varepsilon_j,\delta_j}(v_{j,i},Q(i\delta_j,\rho^{i}_{j}))=r^{d-p}_{\delta_j}\mathcal{F}^{T}_{s_j,\delta_j}(\tilde{v}_{j,i},Q_{R^{i}_{j}})=r^{d-p}_{\delta_j}\varphi_{s_j,T,R^{i}_{j}}(\tilde{u}^{i}_{j}).
        \end{equation}
        
            \noindent Since $\tilde{v}_{j,i}\equiv 0$ on $Q_{1}$, then $u_{j}\equiv 0$ on $Q(i\delta_j,r_{\delta_j})$ for every $i \in Z_j(\Omega')$.
             We need to prove that $u_j \rightarrow u$ in $L^{p}(\Omega;\mathbb{R}^{m})$. By the fact that $\tilde{u}_j\rightarrow u$ in $L^{p}(\Omega';\mathbb{R}^{m})$ and $w_j \rightarrow u$ in $L^{p}(\Omega';\mathbb{R}^{m})$, we have
            \begin{equation*}
                \begin{split}
            \int_{\Omega'}|u_j-w_j|^{p}&=\sum_{\alpha \in \Omega'_{j}}\varepsilon_j^{d}|u_{j}(\alpha)-w_j(\alpha)|^{p}+o(1)\\
                    &\leq \left[\sum_{\alpha \in E_j}\varepsilon_j^{d}|u_{j}(\alpha)-w_j(\alpha)|^{p}\right]+o(1)\\
                    &\leq C\bigg[\int_{\Omega'}|w_{j}-\tilde{u}_{j}|^{p}\,dx+\sum_{i \in Z_{j}(\Omega')}\sum_{\alpha \in Q(i\delta_{j},\rho^{i}_{j})}\varepsilon_{j}^{d}|w_{j}(\alpha)-v_{j,i}(\alpha)|^{p}\bigg]+o(1).
                \end{split}
            \end{equation*}
            For $\alpha \in Q(i\delta_j,\rho^{i}_{j})$ set 
            \begin{equation*}
                g_{j,i}(\alpha):=w_j(\alpha)-v_{j,i}(\alpha).
            \end{equation*}
            Since $v_{j,i}(\alpha)=\tilde{u}^{i}_{j}$ and $w_{j}=\tilde{u}^{i}_{j}$ on $\partial^{\varepsilon_j T}Q(i\delta_j,\rho^{i}_{j})$, it is clear that $g_{i,j}\in W^{1,p}_{0}(Q(i\delta_j,\rho^{i}_{j});\mathbb{R}^{m})$. By \eqref{recovery sequence}, \eqref{minimizing sequence}, assumption (G0), Proposition \ref{lastool} and a discrete version of Poincarè inequality  we get
            \begin{equation*}
            \begin{split}
                 \sum_{i \in Z_{j}(\Omega')}\sum_{\alpha \in Q(i\delta_{j},\rho^{i}_{j})}&\varepsilon_{j}^{d}|g_{j,i}(\alpha)|^{p}\leq C \delta_j^{p}\sum_{i \in Z_j(\Omega')}\sum_{k=1}^{d}\sum_{\alpha \in (Q(i\delta_j,\rho^{i}_{j}))_{\varepsilon_j}(e_k)}\varepsilon_j^{d}\left|\frac{g_{j,i}(\alpha+\varepsilon_je_{k})-g_{j,i}(\alpha)}{\varepsilon_j}\right|^{p}\\
                & \leq C \left[\delta_{j}^{p}\sum_{i \in Z_{j}(\Omega')}G^{1,p}_{\varepsilon_j}(w_{j},Q(i\delta_j,\rho^{i}_{j}))+\delta_j^{p}\sum_{i \in Z_{j}(\Omega')}G^{1,p}_{\varepsilon_{j}}(v_{j,i},Q(i\delta_{j},\rho^{i}_{j}))\right]\\
                &\leq C \left[\delta_{j}^{p}\mathcal{F}^{T}_{\varepsilon_j}(w_{j},\Omega')+\delta_{j}^{p}\sum_{i \in Z_{j}(\Omega')}\mathcal{F}^{T}_{\varepsilon_j,\delta_j}(v_{j,i},Q(i\delta_j,\rho^{i}_{j}))\right]\\
                &\leq C\delta_j^{p}\left[\int_{\Omega'}f^{T}_{\hom}(\nabla u)+\sum_{i \in Z_{j}(\Omega')}r_{\delta_{j}}^{d-p}\varphi_{s_{j},T,R^{i}_{j}}(\tilde{u}^{i}_{j})+\frac{1}{k}\right]\\
                &=C\delta_{j}^{p}\left[\int_{\Omega'}f^{T}_{\hom}(\nabla u)+\gamma^{d-p}\sum_{i \in Z_{j}(\Omega')}\delta_{j}^{d}\varphi_{s_{j},T,R^{i}_{j}}(\tilde{u}^{i}_{j})+\frac{1}{k}\right]\\
                &=o(1)
            \end{split}
            \end{equation*}

\noindent Finally we pass to the estimate of the energy. By (\ref{cambio di variabili sulle perforazioni}) we get 
\begin{equation*}
    \begin{split}
        F^{T}_{\varepsilon_j,\delta_j}(u_j)&\leq \mathcal{F}^{T}_{\varepsilon_j}(w_j,\Omega')+\sum_{i \in Z_{j}(\Omega')}F^{T}_{\varepsilon_j,\delta_j}(v_{j,i},Q(i\delta_j,\rho^{i}_{j}))\\
        &=\mathcal{F}^{T}_{\varepsilon_j}(w_j,\Omega')+\sum_{i \in Z_j(\Omega')}r_{\delta_j}^{d-p}\varphi_{s_j,T,R^{i}_{j}}(\tilde{u}^{i}_{j}).
    \end{split}
\end{equation*}
Applying Proposition \ref{lastool}, (\ref{recovery sequence}) and (\ref{inequality JL}) we get
\begin{equation*}
   \limsup_{j \rightarrow +\infty} F^{T}_{\varepsilon_j,\delta_j}(u_j)\leq \int_{\Omega'}f^{T}_{\hom}(\nabla u)\, dx +\gamma^{d-p}\int_{\Omega'}\varphi^{T}(u)\, dx+\frac{C}{k}.
\end{equation*}
Then, the conclusion follows by the arbitrariness of $k \in \mathbb{N}$ and letting $\Omega' \rightarrow \Omega$.

\end{document}